\documentclass{amsart}

\makeatletter
\def\subsection{\@startsection{subsection}{2}%
  \z@{.5\linespacing\@plus.3\linespacing}{.3\linespacing}%
  {\normalfont\bfseries}}
\makeatother

\makeatletter
\def\subsubsection{\@startsection{subsubsection}{3}%
  \z@{1.2ex\@plus0.8ex}{-1ex}%
  {\normalfont\itshape}}%
\makeatother

\usepackage[defaultlines=4,all]{nowidow} 

\usepackage{xcolor}
\definecolor{dark-blue}{rgb}{0.15,0.15,0.4}

\usepackage[hyphenbreaks]{breakurl} 
\usepackage[colorlinks=true, linkcolor=black, citecolor=black, urlcolor=dark-blue, anchorcolor=dark-blue]{ hyperref }

\usepackage{ mathtools, amssymb, stmaryrd, mathdots, bbm, mleftright, faktor }

\usepackage{ letltxmacro }    
\usepackage{ xparse }
\usepackage{ stackengine, graphicx, scalerel }

\usepackage{ tikz, tikz-cd }
\usetikzlibrary{shapes.misc}

\usepackage{ booktabs }

\usepackage[shortlabels]{ enumitem }

\usepackage{ standalone, graphicx }
\usepackage[section]{ placeins } 
\usepackage[labelformat=simple]{ subcaption }

\usepackage[margin=1cm]{ caption }

\usepackage{ amsthm, thmtools, hyperref, url }
\usepackage[nameinlink]{ cleveref } 

\numberwithin{equation}{section}
\declaretheorem[style=plain,numberlike=equation]{theorem}
\declaretheorem[style=plain,numberlike=theorem]{lemma}
\declaretheorem[style=plain,numberlike=theorem]{proposition}
\declaretheorem[style=plain,numberlike=theorem]{corollary}
\declaretheorem[style=remark,numberlike=theorem]{remark}
\declaretheorem[style=definition,numberlike=theorem]{definition}
\declaretheorem[style=definition,numberlike=theorem]{example}
\declaretheorem[name=Example,style=definition,numbered=no]{example*}
\declaretheorem[style=plain,numberwithin=theorem,name=Claim]{subclaim}
  
\crefformat{section}{\S#2#1#3}
\crefrangeformat{section}{\S#3#1#4--#5#2#6}
\crefmultiformat{section}{\S#2#1#3}{ and~\S#2#1#3}{, \S#2#1#3}{, and~\S#2#1#3}
\crefrangemultiformat{section}{\S#2#1#3}{ and~\S#2#1#3}{, \S#2#1#3}{, and~\S#2#1#3}

\newcommand{\onto}{\twoheadrightarrow}

\NewDocumentCommand\set{s m}{%
    \IfBooleanTF#1%
    {\left\{ #2 \right\}}%
    {\{#2\}}%
}
\NewDocumentCommand\setbuild{s m m}{%
    \IfBooleanTF#1%
    {\ensuremath{\left\{\, #2 \, \middle| \, #3 \,\right\}}}%
    {\ensuremath{\{\, #2 \, \mid \, #3 \,\}}}%
}
\NewDocumentCommand\spangle{s m m m}{%
    \IfBooleanTF#1%
    {\ensuremath{\left\langle\, #2 \, \middle| \, #3 \,\right\rangle_{#4}}}%
    {\ensuremath{\langle\, #2 \, \mid \, #3 \,\rangle_{#4}}}%
}

\newcommand{\symdiff}{\mathbin{\triangle}} 

\DeclarePairedDelimiter{\abs}{\lvert}{\rvert}

\newcommand{\N}{\mathbb{N}}

\newcommand{\iso}{\cong} 
\newcommand{\tensor}{\otimes}

\DeclareMathOperator{\charac}{char}

\DeclareMathOperator{\im}{im}

\newcommand{\restrictto}[1]{\mathord{|}_{#1}}
\newcommand{\blank}{{-}}    

\DeclareMathOperator{\sign}{sgn}

\newcommand{\Y}[1]{[#1]}

\newcommand{\CSYT}{\mathrm{CSYT}}
\newcommand{\SSYT}{\mathrm{SSYT}}

\newcommand{\col}{\mathrm{col}}
\newcommand{\row}{\mathrm{row}}

\newcommand{\rt}[1]{[#1]} 
\newcommand{\act}[1]{|#1|} 
\newcommand{\sct}[1]{||#1||} 
\newcommand{\polyt}[1]{\mathrm{e}(#1)} 

\newcommand{\ppa}{\cdot}        

\newcommand{\sym}{\mathrm{sym}}
\newcommand{\Sk}{\mathsf{Sk}}

\newcommand{\Tbx}[2]{\mathsf{Tbx}^{#1}(#2)}
\newcommand{\Tsym}[2]{\mathsf{Tbx}_{\sym}^{#1}(#2)}

\newcommand{\Gar}{\mathsf{GR}}
\newcommand{\GR}[2]{\Gar^{#1}(#2)}
\newcommand{\GRel}[1]{\mathsf{R}_{#1}}

\newcommand{\symGR}[2]{\Gar_{\sym}^{#1}(#2)}

\newcommand{\skGR}[2]{\Sk\GR{#1}{#2}}
\newcommand{\skR}{\mathsf{R}^{\Sk}}
\newcommand{\skGRel}[1]{\skR_{#1}}

\newcommand{\skewtabs}[2]{\Sk^{#1'} #2}

\newcommand{\alttabs}[2]{\Wedge^{#1'} #2}
\newcommand{\symalttabs}[2]{\Wedge_{\sym}^{#1'} #2}

\newcommand{\polytabspace}[2]{\Schur{#1}(#2)}
\newcommand{\sympolytabspace}[2]{\Schur{#1}_{\sym}(#2)}

\newcommand{\sfun}{\mathcal{F}}
\newcommand{\gtensor}{\mathcal{G}_{\otimes}}

\DeclareMathOperator{\Sym}{Sym}
\newcommand{\wwedge}{\mathchoice{{\textstyle\bigwedge}}%
    {{\bigwedge}}%
    {{\textstyle\wedge}}%
    {{\scriptstyle\wedge}}}
\DeclareMathOperator{\Wedge}{\wwedge}
\newcommand{\Schur}[1]{\nabla^{#1}}

\newcommand{\tless}{\triangleleft}

\DeclareMathOperator{\Hom}{Hom}

\DeclareMathOperator{\soc}{soc}

\newcommand{\Fr}{\mathsf{Fr}}

\renewcommand*{\det}{\qopname\relax o{det}}

\DeclareMathOperator{\GL}{GL}

\newcommand{\twobytwosmallmatrix}[4]{\ensuremath{%
\begin{psmallmatrix}
    #1 & #2 \\ #3 & #4 \\
\end{psmallmatrix} %
}}

\newcommand{\binomII}[2]{
\left.\mathchoice
  {\left(\kern-0.48em\binom{\smash{#1}}{\smash{#2}}\kern-0.48em\right)}
  {\big(\kern-0.30em\binom{\smash{#1}}{\smash{#2}}\kern-0.30em\big)}
  {\left(\kern-0.30em\binom{\smash{#1}}{\smash{#2}}\kern-0.30em\right)}
  {\left(\kern-0.30em\binom{\smash{#1}}{\smash{#2}}\kern-0.30em\right)}
\right.}

\renewcommand{\epsilon}{\varepsilon}
\renewcommand{\phi}{\varphi}
\renewcommand{\leq}{\leqslant}

\renewcommand{\geq}{\geqslant}

\newcommand{\entryset}{\mathcal{B}}
\newcommand{\dWeyl}[1]{\nabla^{#1}(E)}
\newcommand{\natpermrep}{W}
\newcommand{\qmap}{q}

\usepackage{ytableau_mod}
\ytableausetup{centertableaux}
\newcommand{\ytabsmall}[1]{%
\ytableausetup{smalltableaux}%
\ytableaushort{#1}%
\ytableausetup{nosmalltableaux}%
}
\newcommand{\ytab}[1]{%
\ytableaushort{#1}%
}
\newcommand{\yrowtab}[1]{%
\ytableausetup{tabloids}%
\ytableaushort{#1}%
\ytableausetup{notabloids}%
}
\newcommand{\ycoltab}[1]{%
\ytableausetup{coltabloids}%
\ytableaushort{#1}%
\ytableausetup{notabloids}%
}
\newcommand{\ycoltabsmall}[1]{%
\ytableausetup{coltabloids,smalltableaux}%
\ytableaushort{#1}%
\ytableausetup{notabloids,nosmalltableaux}%
}
\newcommand{\yskewcoltab}[1]{%
\ytableausetup{skewcoltabloids}%
\
\ytableaushort{#1}%
\
\ytableausetup{notabloids}%
}
\newcommand{\yskewcoltabsmall}[1]{%
\ytableausetup{skewcoltabloids,smalltableaux}%
\ytableaushort{#1}%
\ytableausetup{notabloids,nosmalltableaux}%
}

\title[The image of the Specht module under the inverse Schur functor]
{The image of the Specht module under the inverse Schur functor in arbitrary characteristic}
\date{}
\author{Eoghan McDowell}

\usepackage[foot]{ amsaddr }

\makeatletter
\def\@setfoot@addresses{%
    \ifx\@affiliation\@empty\else
    \affiliationname\ \@affiliation \@addpunct.%
    \fi
    \ifx\web@page\@empty\else
    
    \webpagename\ \web@page \@addpunct.%
    \fi
}
\def\affiliationname{\textit{Affiliation}:}
\def\affiliation#1{\gdef\@affiliation{#1}}
\let\@affiliation\@empty
\def\webpagename{\textit{Webpage}:}
\def\webpage#1{\gdef\web@page{#1}}
\let\web@page\@empty
\makeatother

\affiliation{Royal Holloway, University of London}
\email{eoghan.mcdowell.2018@rhul.ac.uk}
\webpage{\href{https://eoghanjmcdowell.com}{eoghanjmcdowell.com}}

\makeatletter
\def\@setsubjclass{%
    \vspace{4pt}

    \vspace{-\baselineskip}
    {\itshape\subjclassname.}\enspace\@subjclass\@addpunct.%
}
\makeatother



\makeatletter
\def\@setthanks{%
\vspace{-\baselineskip}\vspace{4pt}
\def\thanks##1{\@par##1\@addpunct.}
\thankses%
}
\def\journalinfo#1{\thanks{#1}}
\makeatother

\journalinfo{
This is the accepted manuscript for an article published in Journal of Algebra 586 (Nov 2021), pp.~865--898, available online
at \href{https://doi.org/10.1016/j.jalgebra.2021.07.013}{doi.org/10.1016/j.jalgebra.2021.07.013}%
}

\begin{document}

\maketitle
\thispagestyle{empty}

\section*{Abstract}

This paper gives a necessary and sufficient condition for the image of the Specht module under the inverse Schur functor to be isomorphic to the dual Weyl module in characteristic \(2\), and gives an elementary proof that this isomorphism holds in all cases in all other characteristics.
These results are new in characteristics \(2\) and \(3\).
We deduce some new examples of indecomposable Specht modules in characteristic \(2\).
When the isomorphism does not hold, the dual Weyl module is still a quotient of the image of the Specht module, and we prove some additional results:
we demonstrate that the image need not have a filtration by dual Weyl modules,
we bound the dimension of the kernel of the quotient map,
and we give some explicit descriptions for particular cases.
Our method is to view the Specht and dual Weyl modules as quotients of suitable exterior powers by the Garnir relations.

\section{Introduction}

The main results of this paper are the descriptions, stated below, of the image of the Specht module under the left-adjoint to the Schur functor.
Throughout, \(n\) and \(d\) denote positive integers, \(\lambda\) denotes a partition of \(n\), and \(K\) denotes a field which may be of characteristic \(0\) or of prime characteristic \(p\).
The Specht module for the symmetric group \(S_n\) labelled by \(\lambda\) is denoted \(S^\lambda\); the dual Weyl module for the general linear group \(\GL_d(K)\) labelled by \(\lambda\) is denoted \(\dWeyl{\lambda}\) (being obtained by applying a certain functorial construction \(\nabla^\lambda\) to the \(d\)-dimensional natural representation \(E\) of \(\GL_d(K)\); see \Cref{subsection:polytabs} for details).
The Schur functor, denoted \(\sfun\), is a functor from the category of polynomial representations of \(\GL_d(K)\) to the category of representations of \(S_r\); its left-adjoint, denoted \(\gtensor\), is right-inverse to \(\sfun\) (see \Cref{section:def_of_Schur_functor} for definitions).

\begin{restatable}{theorem}{mainresultnottwo}\label{thm:mainresultnot2}
Suppose \(K\) has characteristic not \(2\).
Then there is an isomorphism \(\gtensor(S^\lambda) \iso \dWeyl{\lambda}\). 
\end{restatable}

\begin{restatable}{theorem}{mainresulttwo}\label{thm:mainresult2}
Suppose \(K\) has characteristic \(2\).
There is a surjection \(\gtensor(S^\lambda) \onto \dWeyl{\lambda}\), which is an isomorphism if \(\lambda\) is \(2\)-regular, or if \(\lambda_1 = \lambda_2 \geq \lambda_3 + 2\) and \(\lambda\) minus its first part is \(2\)-regular.
Supposing also \(d \geq n-2\), if \(\lambda\) is not of this form then the surjection is not an isomorphism.
\end{restatable}

Although the definition of \(\sfun\) requires \(K\) to be infinite and \(d \geq n\), and the term ``dual Weyl module'' is non-standard when \(K\) is finite, the definitions of \(\gtensor\) and of \(\nabla^\lambda\) do not have these restrictions, and neither do our main theorems.

The isomorphism \(\gtensor(S^\lambda) \iso \dWeyl{\lambda}\) is known to hold in all characteristics other than \(2\) and \(3\) (\cite[3.2]{Kleshchev2001cohomology}),
and more generally in the context of \(q\)-Schur algebras and Hecke algebras of quantum characteristic at least \(4\) (\cite[Theorem 3.4.2]{HemmerNakano2004hecke}).
A related result identifying the image of the twisted Young module in characteristics other than \(2\) is given in \cite[Theorem 5.2.4]{CPS1996stratifying}.

The novelty of these results, therefore, is in characteristics \(2\) and \(3\).
Additionally, our approach establishes the isomorphisms in characteristics other than \(2\) with much less machinery than the accounts cited above.

The viewpoint taken in this paper yields a concrete model for the image \(\gtensor(S^\lambda)\) as a quotient of a skew symmetric power by a modified set of Garnir relations (see \Cref{lemma:gtensor_applied_to_sym_ses}).
However, it can be challenging to understand this module further when it is not isomorphic to \(\nabla^\lambda(E)\) (necessarily in characteristic \(2\)).
We show that it need not have a filtration by dual Weyl modules (\Cref{eg:no_nabla_filtration}).
Even when \(d=1\), the description (\Cref{prop:when_d=1}) is nontrivial: unlike the dual Weyl module, the image \(\gtensor(S^\lambda)\) can be nonzero for partitions of arbitrary length.
Our other results include
bounding the dimension of the kernel of the map in \Cref{thm:mainresult2} as \(O(d^{n-1})\) for fixed \(n\) as \(d\) varies (\Cref{prop:dim_growth_of_U});
describing \(\gtensor(S^\lambda)\) when \(\lambda\) is a hook partition and \(d=2\) (\Cref{prop:hooks_when_d=2});
and identifying the composition factors of \(\gtensor(S^\lambda)\) for partitions of \(n \leq 5\) (\Cref{eg:comp_factors_for_partitions_of_n_less_than_5}).

We deduce from \Cref{thm:mainresult2} the following corollary on the indecomposability of some Specht modules.
In characteristics other than \(2\), all Specht modules are known to be indecomposable, and in characteristic \(2\) those indexed by \(2\)-regular partitions are known to be indecomposable.
Determining the decomposability of the remaining Specht modules is a difficult open problem.
Families of decomposable Specht modules have been identified by Murphy \cite{Murphy1980decomposableSpechts}, Dodge and Fayers \cite{Dodge2012decomposableSpechts}, and Donkin and Geranios \cite{DONKIN2020decompositions}.

\begin{restatable}{corollary}{indecSpechts}\label{cor:indec_Spechts}
Suppose \(K\) is infinite and has characteristic \(2\).
Let \(\lambda\) be a partition such that \(\lambda_1 = \lambda_2 \geq \lambda_3 +2\) and such that \(\lambda\) minus its first part is \(2\)-regular.
Then \(S^\lambda\) is indecomposable. 
\end{restatable}


Note that the Schur functor also has a right-adjoint right-inverse, denoted \(\mathcal{G}_{\Hom}\), which is related to the left-adjoint by duality (see \cite[2.2]{COHEN20101419}).
We consider only the functor \(\gtensor\); results for \(\mathcal{G}_{\Hom}\) can be obtained from ours via this duality. 

This paper is structured as follows.
In \Cref{section:defs_and_cons} we describe \(S^\lambda\) and \(\nabla^\lambda(E)\) as spaces of polytabloids and as quotients by the Garnir relations.
In \Cref{section:def_of_Schur_functor} we define the Schur functor and its inverse.
In \Cref{section:quotient_construction} we introduce skew column tabloids and skew Garnir relations and use them to model \(\gtensor(S^\lambda)\), obtaining \Cref{thm:mainresultnot2}. 
In \Cref{section:skew_Garnir_relations} we prove a key technical result on the skew Garnir relations in characteristic \(2\).
In \Cref{section:image_of_Specht_in_characteristic_2} we deduce \Cref{thm:mainresult2} and prove the additional results in characteristic \(2\) mentioned above.

\section{Polytabloid constructions for the Specht and dual Weyl modules}
\label{section:defs_and_cons}

In this section we construct the Specht and dual Weyl modules both as modules consisting of \emph{polytabloids} (that is, column-antisymmetrised linear combinations of row tabloids), and as quotients by the \emph{Garnir relations}.
This is an adaptation of James's construction of the Specht modules \cite{gdjames1978reptheorysymgroups} to an arbitrary group \(G\), given for the general linear group in \cite{boeck2018plethysms}.

We perform this construction using tableaux whose entries are elements of an ordered basis \(\entryset\) for a left \(KG\)-module \(V\).
The \(G\)-action on \(V\) induces a ``diagonal'' left \(G\)-action on the space of tableaux (and their equivalence classes) by entrywise action and multilinear expansion; the \(K\)-vector spaces defined in this section thus become left \(KG\)-modules.
We denote these group actions by concatenation.

For our purposes, either \(G = S_n\) and \(V\) is the natural permutation module \(\natpermrep\), or \(G = \GL_d(K)\) and \(V\) is the natural \(d\)-dimensional module \(E\).
Nevertheless, the constructions in this section (and the results of \Cref{section:skew_Garnir_relations} later) are valid for any choice of group and module.

\subsection{Tableaux and tabloids}
\label{subsection:tableaux_and_tabloids}

In this subsection we define the combinatorial objects necessary for the construction of our representations.

\subsubsection{Partitions}

A \emph{partition} of \(n\) is a weakly decreasing sequence of positive integers whose sum is \(n\).
If \(\lambda\) is a partition with \(l\) parts, we interpret \(\lambda_i = 0\) for \(i > l\).
The \emph{conjugate} (or \emph{transpose}) of \(\lambda\), denoted \(\lambda'\), is the partition defined by \(\lambda'_i = \abs{\setbuild{j \geq 1}{ \lambda_j \geq i}}\).

The \emph{Young diagram} of \(\lambda\) is the set \(\Y{\lambda} = \setbuild{(i,j)}{1 \leq i \leq \lambda'_1, 1 \leq j \leq \lambda_i}\), which we picture lying in the plane (with the \(x\)-direction downward and the \(y\)-direction rightward).
An element of a Young diagram is called a \emph{box}.
Let \(\row_{i}[\lambda]\) and \(\col_{j}{[\lambda]}\) denote the sets of boxes in row \(i\) and column \(j\) of \(\Y{\lambda}\) respectively.

\subsubsection{Tableaux}

A \emph{tableau} of shape \(\lambda\) with entries in \(\entryset\) is a function \(\Y{\lambda} \to \entryset\).
The image of a box \(b \in \Y{\lambda}\) under a tableau \(t\) is the \emph{entry} of \(t\) in \(b\).
We depict a tableau \(t\) by filling the boxes in the Young diagram of \(\lambda\) with their entries in \(t\). 
The \emph{weight} of \(t\) is the multiset of entries of a tableau \(t\), expressed as a composition of \(n\) via the total ordering on \(\entryset\).

Let \(\Tbx{\lambda}{V}\) be the \(KG\)-module with basis the set of tableaux of shape \(\lambda\) with entries in \(\entryset\).
There is a (non-unique) isomorphism \(\Tbx{\lambda}{V} \iso V^{\tensor n}\).

If the entries of a tableau strictly increase along the rows or columns, we say it is \emph{row standard} or \emph{column standard} respectively.
If the entries of a tableau weakly increase along the rows or columns, we say it is \emph{row semistandard} or \emph{column semistandard} respectively.

If a tableau is both row standard and column standard, we say it is \emph{standard}.
If a tableau is both row semistandard and column standard, we say it is \emph{semistandard};
the set of semistandard tableaux of shape \(\lambda\) with entries in \(\entryset\) is denoted \(\SSYT_\entryset(\lambda)\).

If a tableau is both row semistandard and column semistandard (known in some other contexts, after reversing the ordering, as a reverse plane partition), we abbreviate this description to \emph{row-and-column semistandard}.

\subsubsection{Tableaux of symmetric type}
\label{subsection:tableaux_of_symmetric_type}

We say a tableau \(t\) is \emph{of symmetric type} when all entries of \(t\) are distinct.
Let \(\Tsym{\lambda}{V}\) be the \(K\)-subspace of \(\Tbx{\lambda}{V}\) spanned by tableaux of symmetric type.
Likewise, for all constructions of spaces in this section, let \(\blank_{\sym}\) denote the construction restricted to tableaux of symmetric type.

Note that these \(K\)-subspaces may not be \(KG\)-submodules in general.
However, if \(V\) is a permutation \(KG\)-module (as in our specialisation to the symmetric group and the natural permutation module), then indeed they are \(KG\)-submodules.

\subsubsection{Place permutation action on tableaux}

Given a set \(X\), we let \(S_X\) denote the group of permutations of \(X\), with permutations written on the right of their arguments.
The group \(S_{\Y{\lambda}}\) of permutations of the Young diagram acts on tableaux on the right by permuting the boxes via
\[
(t \ppa \sigma)(b) = t(b\sigma^{-1})
\]
for \(\sigma \in S_{\Y{\lambda}}\) and \(b \in \Y{\lambda}\).
This action makes \(\Tbx{\lambda}{V}\) into a \(KS_{\Y{\lambda}}\)-module (moreover, a permutation module).

This action is useful for defining more complicated structures, but \(S_{\Y{\lambda}}\) is (in general) not the group whose representation theory we are interested in.

\newcommand{\RPP}{\mathrm{RPP}}
\newcommand{\CPP}{\mathrm{CPP}}

Define the sets of \emph{row-preserving} and \emph{column-preserving place permutations}, subgroups of \(S_{\Y{\lambda}}\), by
\[
    \RPP(\lambda) = \prod_{i=1}^{\lambda'_1} S_{\row_{i}[\lambda]}
    \qquad\quad \text{and} \qquad\quad
    \CPP(\lambda) = \prod_{j=1}^{\lambda_1} S_{\col_{j}[\lambda]}.
\]

\subsubsection{Row tabloids}

\newcommand{\Jsym}{J_{\mathsf{Sym}}}

A \emph{row tabloid} is an equivalence class of tableaux under row equivalence.
Concretely, we quotient the space of tableaux \(\Tbx{\lambda}{V}\) by the subspace
\begin{align*}
    \Jsym &= \spangle{x \ppa \sigma - x}{x \in \Tbx{\lambda}{V},\, \sigma \in \RPP(\lambda)}{K},
\end{align*}
and say the row tabloid corresponding to a tableau \(t\) is the element \(t + \Jsym\) in the quotient \(\Tbx{\lambda}{V}/\Jsym\).
We write the row tabloid corresponding to \(t\) as \(\rt{t}\), and draw a row tabloid \(\rt{t}\) by deleting the vertical lines from a drawing of \(t\), as depicted below in the case \(\lambda=(3,2)\).
\[
t = \ytab{124,35} 
\quad\implies\quad
\rt{t} = \yrowtab{124,35}
\]

By construction, \(\rt{t \ppa \sigma} = \rt{t}\) for any \(\sigma \in \RPP(\lambda)\).
Moreover the space of row tabloids is naturally isomorphic as a \(KG\)-module to the symmetric power \(\Sym^\lambda V = \bigotimes_{i=1}^{\lambda'_1} \Sym^{\lambda_i} V\), where ``symmetric power'' refers to the quotient of the tensor power
\[
\Sym^r V \iso \faktor{V^{\tensor r}}{\spangle{
            w \cdot \sigma - w
        }{
            w \in V^{\otimes r},\, \sigma \in S_r
        }{K}
    }.
\]
We therefore use \(\Sym^\lambda V\) to denote the space of row tabloids.

\subsubsection{Column tabloids}

\newcommand{\Jsk}{J_{\Sk}}
\newcommand{\Jalt}{J_{\mathsf{Alt}}}

When defining column tabloids, we wish to also associate signs to the equivalence classes.
This is achieved by quotienting the space of tableaux \(\Tbx{\lambda}{V}\) by the subspace
\begin{align*}
    \Jalt &= \spangle{x \in \Tbx{\lambda}{V}}{\text{\(x \ppa \tau = x\) for some transposition \(\tau \in \CPP(\lambda)\)}}{K}.
\end{align*}

The \emph{alternating column tabloid} corresponding to a tableau \(t\) is the element \(t + \Jalt\) in the quotient \(\Tbx{\lambda}{V}/\Jalt\).
We write this tabloid as \(\act{t}\), and draw an alternating column tabloid by deleting the horizontal lines from a drawing of the corresponding tableau, as depicted below in the case \(\lambda = (3,2)\).
\[
    t = \ytab{124,35} 
        \quad \implies\quad
    \act{t} = \ycoltab{124,35}
\]

Observe that \(\act{t \ppa \sigma} = \act{t} \sign(\sigma)\) for any \(\sigma \in \CPP(\lambda)\),
and furthermore \(\act{t} = 0\) if \(t\) has a repeated entry in a column. 
For example, with \(\lambda = (1,1)\), the elements \(\,\ytabsmall{1,1}\,\) and \(\,\ytabsmall{1,2} + \ytabsmall{2,1}\,\) of \(\Tbx{\lambda}{V}\) are fixed by the transposition swapping the only two boxes, so these element lies in \(\Jalt\) and hence \(\ycoltabsmall{1,1} = 0\) and \(\ycoltabsmall{1,2} = -\,\ycoltabsmall{2,1}\).
(To see that \(\act{t \ppa \sigma} = \act{t} \sign(\sigma)\) when \(\sigma \in \CPP(\lambda)\) is a product of several transpositions, consider the collection of elements of the form \(t \ppa \tau_1 \cdots \tau_i \, + \, t \ppa \tau_1 \cdots \tau_{i-1} \in \Jalt\) where \(\tau_1, \tau_2, \ldots\) is a sequence of transpositions whose product is \(\sigma\).)

The space of alternating column tabloids is therefore naturally isomorphic as a \(KG\)-module to the exterior power \(\alttabs{\lambda}{V} = \bigotimes_{i=1}^{\lambda_1} \Wedge^{\lambda'_i}(V)\),
and we use \(\alttabs{\lambda} V\) to denote the space of alternating column tabloids.
This space has \(K\)-basis \(\setbuild{\act{t}}{t \in \CSYT(\lambda)}\).

When we come to model the image of the Specht module under the inverse Schur functor, we introduce a different form of column tabloid called a \emph{skew column tabloid} (\Cref{def:skew_column_tabloid}).

\subsubsection{Column ordering on tableaux}
\label{subsection:orders_on_tableaux}

\newcommand{\colless}{<_{\mathrm{c}}}
\newcommand{\colsim}{\sim_{\mathrm{c}}}
\newcommand\scalesim{\scalebox{.8}{$\SavedStyle\sim$}}
\newcommand{\colleq}{%
    \mathrel{\ensurestackMath{\ThisStyle{%
        \stackengine{-.6\LMex}{\SavedStyle{<_{\mathrm{c}}}}{%
            \rotatebox{-25}{%
                \scalesim
            }%
        }{U}{l}{F}{T}{S}
    }}}%
}

We make use of an ordering on tableaux which we call the column ordering \(\colless\).
When considering only tableaux of symmetric type, \(\colless\) is the column analogue of the order defined by James in \cite[Definition~3.10]{gdjames1978reptheorysymgroups} (though ours is defined on tableaux rather than tabloids).
The order is much easier to interpret in this case:
to compare two distinct tableaux of symmetric type, identify the largest entry which does not appear in the same column in both tableaux, and declare the \(\colless\)-greater tableau to be the one for which this element is further left.
We illustrate the \(\colless\)-least and \(\colless\)-greatest standard tableaux of symmetric type in the case \(\lambda = (4^3,2,1)\) and \(\entryset = [15]\) in \Cref{fig:col-extremal_tableaux}.
It can be shown that \(\colless\) extends the dominance order \(\tless\) defined on column equivalence classes in \cite[Definition~13.8]{gdjames1978reptheorysymgroups}, in the sense that \(t/{\colsim} \tless u/{\colsim}\) implies \(t \colless u\) (where \(t/{\colsim}\) denotes the column equivalence class of \(t\)).

\begin{definition}
\label{def:column_ordering}
The \emph{column ordering}, \(\colless\), is the strict partial order on the set of tableaux of a fixed shape defined as follows.
Consider tableaux \(t\) and \(u\) of shape \(\lambda\).
\begin{itemize}
    \item
If there is equality \(\col_{j}(t) = \col_{j}(u)\) (as multisets) for all \(1 \leq j \leq \lambda_1\), then \(t\) and \(u\) are \(\colless\)-incomparable, and we write \(t \colsim u\).
    \item 
Otherwise, let \(m \in \entryset\) be maximal such that there exists \(j\) such that \(m \in \col_{j}(t) \symdiff \col_{j}(u)\) (where \(\symdiff\) denotes the multiset symmetric difference).
Let \(j\) be minimal such that \(m \in \col_{j}(t) \symdiff \col_{j}(u)\).
If \(m \in \col_{j}(t) \setminus \col_{j}(u)\), then \(u \colless t\); conversely
if \(m \in \col_{j}(u) \setminus \col_{j}(t)\), then \(t \colless u\).
\end{itemize}
We write \(t \colleq u\) to mean \(t \colless u\) or \(t \colsim u\).
\end{definition}

The relation \(\colsim\) is an equivalence relation.
Since tableaux \(t\) and \(u\) are \(\colless\)-incomparable if and only if \(t \colsim u\), the relation \(\colleq\) is a total preorder, also known as a weak order (that is, \(\colleq\) is a partial order with the antisymmetry requirement relaxed, permitting \(t \colleq u\) and \(u \colleq t\) to hold simultaneously for distinct \(t\) and \(u\), and with the property that at least one of \(t \colleq u\) and \(u \colleq t\) holds for any pair of tableaux \(t\) and \(u\)).

\begin{figure}[ht]
    \begin{subfigure}{0.49\linewidth}
\centering\ytableausetup{nosmalltableaux}
\ytableaushort{16{10}{13},27{11}{14},38{12}{15},49,5}
\caption{\(\colless\)-least.}
    \end{subfigure}
    \begin{subfigure}{0.49\linewidth}
\centering
\ytableaushort{1234,5678,9{10}{11}{12},{13}{14},{15}}
\caption{\(\colless\)-greatest.}
    \end{subfigure}
    \caption{Extremal standard tableaux for \(\lambda = (4^3,2,1)\) and \(\entryset = [15]\).}
    \label{fig:col-extremal_tableaux}
\end{figure}

\subsection{Polytabloids and the Specht and dual Weyl modules}
\label{subsection:polytabs}

The \emph{polytabloid} corresponding to a tableau \(t\) is the element of \(\Sym^{\lambda}{V}\) given by
\begin{align*}
    \polyt{t} &= \sum_{\sigma \in \CPP(\lambda)} \rt{t \cdot \sigma} \sign{\sigma}.
\end{align*}

The \(K\)-subspace, and moreover \(KG\)-submodule, of \(\Sym^{\lambda}{V}\) generated by the polytabloids is denoted \(\polytabspace{\lambda}{V}\).
Since the action of \(G\) commutes with the place permutation action, it is straightforward to verify that \(\polytabspace{\lambda}{V}\) is a \(KG\)-submodule of \(\Sym^{\lambda}{V}\).
This space has the following well-known basis.

\begin{proposition}[{\cite[Proposition 2.11]{boeck2018plethysms}}]
\label{prop:basis_for_polytabs}
The set \(\setbuild{\polyt{s}}{s \in \SSYT(\lambda)}\) is a \(K\)-basis for \(\polytabspace{\lambda}{V}\).
\end{proposition}

When \(G = S_n\) and \(V = \natpermrep\) is the \(n\)-dimensional natural permutation module for \(S_n\), the \(KS_n\)-module \(\sympolytabspace{\lambda}{\natpermrep}\) of polytabloids of symmetric type is known as a \emph{Specht module} and is denoted \(S^\lambda\).

When \(G = \GL_d(K)\) and \(V = E\) is the \(d\)-dimensional natural \(K\GL_d(K)\)-module, we call the \(KG\)-module \(\polytabspace{\lambda}{E}\) of polytabloids a \emph{dual Weyl module}.
Green \cite[Section 4]{green2006polynomial} denotes this module \(D_{\lambda,K}\), and constructs it as a module for the Schur algebra spanned by bideterminants; it is the dual of what he calls the Weyl module \cite[Section 5]{green2006polynomial}.
James \cite[Definitions 17.2, 17.4 and 26.4, pp.~65,127,129]{gdjames1978reptheorysymgroups} denotes this module \(W^\lambda\), and constructs it by summing the images of the space of polytabloids of symmetric type under maps which induce each possible weight.

\begin{remark}
The construction of the space of polytabloids defines an endofunctor \(\nabla^\lambda\) on the category of \(KG\)-modules for any group \(G\).
These endofunctors are sometimes known as \emph{Schur functors}, but are not to be confused with the Schur functor \(\sfun\) defined in \Cref{section:def_of_Schur_functor}.
\end{remark}

\subsection{Garnir relations and quotient construction of Specht and dual Weyl modules}
\label{section:basis_for_Garnir_relations}

An immediate consequence of the definition of a polytabloid is that \( \polyt{t \cdot \sigma} = \polyt{t} \sign{\sigma}\) for \(\sigma \in \CPP(\lambda)\),
and that \(\polyt{t} = 0\) if \(t\) has a repeated entry in a column.
It follows that the map \(e \colon \alttabs{\lambda}{V} \to \polytabspace{\lambda}{V}\) defined by \(K\)-linear extension of
\begin{align*}
    e \colon \act{t} &\mapsto \polyt{t}
\end{align*}
is well-defined and surjective.
It is also \(G\)-equivariant.
The aim of this subsection is to describe the kernel of \(e\) explicitly, thus constructing \(\nabla^\lambda V\) as a quotient of \(\Wedge^{\lambda'} V\).
We exhibit the well-known argument that the kernel consists of elements called \emph{Garnir relations}, and furthermore identify a basis for this space.

Our definition of Garnir relations is as certain linear combinations of alternating column tabloids (that is, as certain elements of \(\alttabs{\lambda}{V}\)).
Garnir elements were defined in the context of tableaux of symmetric type by James \cite[Section 7]{gdjames1978reptheorysymgroups} as certain elements of the group algebra \(KS_n\); in that context, these elements yield our notion of a Garnir relation when they act on suitable column tabloids of symmetric type.
The relations used by de Boeck, Paget and Wildon \cite[Lemma 2.4 and Equation 2.5]{boeck2018plethysms} are images of our Garnir relations under the map \(e\).
Fulton \cite[Section 8]{fulton1997youngtableaux} describes a similar collection of linear combinations of alternating column tabloids which he calls quadratic relations; these generate the same \(K\)-subspace of \(\alttabs{\lambda}{V}\) as our Garnir relations.

\begin{definition}[Garnir relations]
\label{def:Garnir_relation}
Let \(t\) be a tableau of shape \(\lambda\) with entries in \(\entryset\).
Let \(1 \leq j < j' \leq \lambda_1\), and let \(A \subseteq \col_j(\lambda)\) and \(B \subseteq \col_{j'}(\lambda)\) be such that \(\abs{A} + \abs{B} > \lambda'_j\).
Choose \(\mathcal{S}\) a set of left coset representatives for \(S_{A} \times S_B\) in \(S_{A \sqcup B}\).
The \emph{Garnir relation} labelled by \((t, A, B)\) is
\[
    \mathsf{R}_{(t, A, B)} = \sum_{\tau \in \mathcal{S}} \act{t \ppa \tau} \sign{\tau}.
\]
Let \(\GR{\lambda}{V}\) denote the subspace, and moreover \(KG\)-submodule, of \(\alttabs{\lambda}{V}\) which is spanned by the Garnir relations.
\end{definition}


In \cite[Lemma 8.4]{gdjames1978reptheorysymgroups}, James shows that his Garnir elements annihilate suitable polytabloids.
The analogue of this result in our approach is the following.

\begin{proposition}[{\cite[Lemma 2.4 and Equation 2.5]{boeck2018plethysms}}]
\label{prop:Garnir_relations_evaluate_to_zero}
Let \(\mathsf{R}_{(t,A,B)}\) be any Garnir relation.
Then \(e(\mathsf{R}_{(t,A,B)}) = 0\).
\end{proposition}

We strengthen this result in \Cref{prop:basis_for_Garnir_relations}, showing that \(\ker e = \GR{\lambda}{V}\).
In the process, we identify a basis for \(\GR{\lambda}{V}\).

As is usual, we generally need only consider Garnir relations in which the chosen columns are adjacent and boxes are taken from the bottom of the left-hand column and the top of the right-hand column, with a single row containing chosen boxes from both columns.
Following the terminology introduced in \cite[Equation 2.5]{boeck2018plethysms}, we call such relations \emph{snake relations} due to the shape of the outline of the chosen boxes.
Formally they are as defined as follows.

\begin{definition}[Snake relations]
\label{def:snake_relation}
A Garnir relation \(\mathsf{R}_{(t,A,B)}\) is called a \emph{snake relation} when, in the notation of \Cref{def:Garnir_relation}, \(j' = j+1\) and there exists \(i\) such that \(A = \setbuild{(r,j)}{i \leq r \leq \lambda'_j}\) and \(B = \setbuild{(r,j')}{1 \leq r \leq i}\).
In this case, we may also label the Garnir relation by \((t, i, j)\).
\end{definition}

\begin{lemma}
\label{lemma:ordered_expression_for_snake_relation}
Let \(t\) be a column standard tableau, and suppose \((i,j)\) is such that \(t(i, j) > t(i, j+1)\).
Then
\[
    \mathsf{R}_{(t,i,j)} = \act{t} + \sum_{\substack{u \colless t}} m_u \act{u}
\]
for some elements \(m_u\) in the subring of \(K\) generated by \(1\).
\end{lemma}

\begin{proof}
By assumption, the sets \(A = \setbuild{(r,j)}{i \leq r \leq \lambda'_j}\) and \(B = \setbuild{(r,j+1)}{1 \leq r \leq i}\) defining the Garnir relation satisfy
\[
    t(1,j+1) < \ldots < t(i,j+1) < t(i,j) < t(i+1, j) < \ldots < t(\lambda'_j, j).
\]
Thus for any \(\sigma \in S_{A \sqcup B}\),
we have \(t \ppa \sigma \colleq t\), with \(t \ppa \sigma \colsim t\) if and only if \(\sigma \in S_{A} \times S_B\).
\end{proof}

Our basis is the following subset of the snake relations.

\begin{definition}[Basic snake relations]
\label{def:basic_snake_relation}
Let \(\Phi\) be a function on column standard tableaux which are not row semistandard whose output on such a tableau \(t\) is a box \((i,j)\) such that \(t(i,j) > t(i,j+1)\).
A snake relation \(\mathsf{R}_{(t,i,j)}\) is called \emph{\(\Phi\)-basic} if \(t\) is column standard but not row semistandard and \((i,j) = \Phi(t)\).
\end{definition}

The purpose of \(\Phi\) is to associate a unique snake relation to each column standard tableau which is not row semistandard.
Any such function suffices: except in the proofs of \Cref{prop:spanning_set_for_skew_Garnir_relations,prop:hooks_when_d=2}, the choice of \(\Phi\) is irrelevant (that is, all the claims, including the statements of those propositions, hold for any choice of \(\Phi\)).
Accordingly, \(\Phi\) is suppressed in the notation.
An example of a suitable function \(\Phi\) is to let \(\Phi(t) = (i,j)\) where \(j\) is least (primarily) and \(i\) is greatest (secondarily) such that \(t(i,j) > t(i,j+1)\); outside the specified proofs, we may consider this to be the function in the definition of basic snake relations.

The proof of the following lemma is essentially the proof of Corollary~2.6 in \cite{boeck2018plethysms}, except that we insist that the chosen snake relations are basic.

\begin{lemma}[{cf.~\cite[Corollary 2.6]{boeck2018plethysms}}]
\label{lemma:standard_expression_for_col_tabloid}
Let \(t\) be any tableau.
Then there exists some linear combination \(\gamma\) of basic snake relations such that
\[
    \act{t} + \gamma = \sum_{s \in \SSYT(\lambda)} a_s \act{s}
\]
for some elements \(a_s\) in the subring of \(K\) generated by \(1\).
\end{lemma}

\begin{proof}
Without loss of generality, we may assume \(t\) is column standard.
If \(t\) is also row semistandard, we are done.
Otherwise, let \((i,j) = \Phi(t)\), and we have that \(t(i,j) > t(i,j+1)\) and that \(\mathsf{R}_{(t,i,j)}\) is a basic snake relation.
By \Cref{lemma:ordered_expression_for_snake_relation}, \(\mathsf{R}_{(t,i,j)} = \act{t} + \sum_{u \colless t} m_u \act{u}\) for some elements \(m_u\) in the subring of \(K\) generated by \(1\).
Then \(\act{t} - \mathsf{R}_{(t,i,j)}\) is a linear combination of column tabloids whose tableaux precede \(t\) in the column ordering.
The result follows by induction.
\end{proof}

\begin{proposition}
\label{prop:basis_for_Garnir_relations}
The basic snake relations form a basis of \(\GR{\lambda}{V}\), and \(\ker e = \GR{\lambda}{V}\).
\end{proposition}

\begin{proof}
From \Cref{prop:Garnir_relations_evaluate_to_zero}, we have that \(\GR{\lambda}{V} \subseteq \ker e\).
That the basic snake relations are \(K\)-linearly independent follows immediately from \Cref{lemma:ordered_expression_for_snake_relation}.
It now suffices to show that the basic snake relations span \(\ker e\).

Let \(\kappa \in \ker e\).
By \Cref{lemma:standard_expression_for_col_tabloid} there exists a linear combination \(\gamma\) of basic snake relations such that
\begin{align*}
    \kappa + \gamma &= \sum_{s \in \SSYT(\lambda)} a_s \act{s}
\intertext{for some elements \(a_s\).
Applying \(e\) to this equation
, we find}
    0 &= \sum_{s \in \SSYT(\lambda)} a_s \polyt{s}.
\end{align*}
But the semistandard polytabloids are \(K\)-linearly independent by \Cref{prop:basis_for_polytabs}, so this implies that \(a_s = 0\) for all \(s\).
Hence \(\kappa = - \gamma\) is in the span of the basic snake relations, as required.
\end{proof}

\begin{corollary}
\label{cor:short_exact_sequence_for_nabla}
There is a short exact sequence
\begin{center}
\begin{tikzcd}
0 \arrow[r] & \GR{\lambda}{V} \arrow[r] & \alttabs{\lambda}{V} \arrow[r, "e"] & \polytabspace{\lambda}{V} \arrow[r] & 0
\end{tikzcd}
\end{center}
in the category of \(KG\)-modules.
\end{corollary}

\subsubsection*{Garnir relations of symmetric type}

We record here that all the results of this section hold upon restriction to the symmetric type subspace.

\begin{definition}[Garnir relations of symmetric type]
We say a Garnir relation labelled by \((t,A,B)\) is \emph{of symmetric type} if \(t\) is of symmetric type.
Let \(\symGR{\lambda}{V}\) denote the subspace of \(\GR{\lambda}{V}\) spanned by the Garnir relations of symmetric type.
\end{definition}

Since the status of being of symmetric type is preserved under the place permutation action, all the summands of a Garnir relation are of symmetric type if and only if the labelling tableau is.
Thus \(\symGR{\lambda}{V}\) is the intersection \(\GR{\lambda}{V} \cap \symalttabs{\lambda}{V}\).
As indicated in \Cref{subsection:tableaux_of_symmetric_type}, the \(K\)-vector spaces \(\symGR{\lambda}{V}\) and \(\symalttabs{\lambda}{V}\) are \(KG\)-modules if \(V\) is a permutation \(KG\)-module, but in general may not be.

By restricting to tableaux of symmetric type, we obtain the following results from the results of the previous subsection.

\begin{proposition}
\label{prop:basis_for_Garnir_relations_of_sym_type}
The basic snake relations of symmetric type form a basis of \(\symGR{\lambda}{V}\).
\end{proposition}

\begin{proposition}
\label{prop:short_exact_sequence_of_symmetric_type}
Suppose \(V\) is a permutation \(KG\)-module and \(\entryset\) is a permutation basis.
Then there is a short exact sequence
\begin{center}
\begin{tikzcd}
0 \arrow[r] & \symGR{\lambda}{V} \arrow[r] & \symalttabs{\lambda}{V} \arrow[r, "e|_{\sym}"] & \sympolytabspace{\lambda}{V} \arrow[r] & 0
\end{tikzcd}
\end{center}
in the category of \(KG\)-modules.
\end{proposition}

\section{The Schur functor and the inverse Schur functor}
\label{section:def_of_Schur_functor}

In this section we introduce the Schur functor and its one-sided inverse.
These functors are described by Green \cite[Section 6]{green2006polynomial} in the language of the Schur algebra; here we give more elementary constructions.

Throughout, \(E\) denotes the natural \(K\GL_d(K)\)-module, of dimension \(d\), and \(\entryset\) denotes its canonical basis. 
We let \(L^\lambda(E) = \soc \dWeyl{\lambda}\) denote the simple \(K\GL_d(K)\)-module indexed by \(\lambda\).
Note that \(\dWeyl{\lambda}\) vanishes if and only if \(\lambda\) has strictly fewer than \(d\) nonzero parts; in this case we interpret \(L^\lambda(E) = 0\).

\begin{definition}[Schur functor]
\label{def:Schur_functor}
Suppose \(K\) is infinite and \(d \geq n\).
The \emph{Schur functor} \(\sfun\) is the functor from the category of polynomial left \(K \GL_d(K)\)-modules of degree \(n\) to the category of left \(K S_n\)-modules.
It is defined by
\[
    \sfun(V) = V_{(1^n,\,0^{d-n})}
\]
where \(V_\alpha\) denotes the \(\alpha\)-weight space of a \(K \GL_d(K)\)-module \(V\).
\end{definition}

\begin{proposition}[{\cite[(6.2a) and (6.3)]{green2006polynomial}}]
\label{prop:properties_of_sfun}
Suppose \(K\) is infinite and \(d \geq n\).
\begin{enumerate}[(i)]
    \item
The functor \(\sfun\) is exact.
    \item
There is an isomorphism \(\sfun(\dWeyl{\lambda}) \iso S^\lambda\).
\end{enumerate}

\end{proposition}

\begin{definition}[Inverse Schur functor]
\label{def:tensor_inverse_Schur_functor}
The \emph{inverse Schur functor} \(\gtensor^{d}\) is a functor from the category of left \(KS_n\)-modules to the category of left \(K\GL_d(K)\)-modules of degree \(n\).
It is defined by
\[
    \gtensor^d(U) = E^{\tensor n} \tensor_{K S_n} U
\]
where \(E^{\tensor n}\) is viewed as a right \(KS_n\)-module via the place permutation action.
We suppress the dependence of \(\gtensor^d\) on \(d\) except where there is need to emphasise it.
Note that unlike the definition of \(\sfun\) we do not require \(K\) to be infinite or \(d \geq n\).
\end{definition}

\begin{proposition}
\label{prop:properties_of_gtensor}
The functor \(\gtensor\) is:
\begin{enumerate}[(i)]
    \item\label{item:gtensor_right_inverse}
right-inverse to \(\sfun\) (that is, \(\sfun\gtensor(V) \iso V\)), provided \(K\) is infinite and \(d \geq n\);
    \item\label{item:gtensor_left_adjoint}
left-adjoint to \(\sfun\), provided \(K\) is infinite and \(d \geq n\);
    \item\label{item:gtensor_right_exact}
right exact.
\end{enumerate}
\end{proposition}

\begin{proof}
For part \ref{item:gtensor_right_inverse}, see \cite[(6.2d)]{green2006polynomial}.
Part \ref{item:gtensor_left_adjoint} is a particular case of the tensor-hom adjunction.
Part \ref{item:gtensor_right_exact} is a general property of tensor functors.
\end{proof}

\newcommand{\allmat}{\operatorname{Mat}}

\subsubsection*{Dimension reduction functor}

We can consider the effect of varying the parameter \(d\) using the following functor, defined in the language of the Schur algebra by Green in \cite[Section 6.5]{green2006polynomial}.
Here we multiply by an idempotent in \(\allmat_d(K)\), the algebra of all (not necessarily invertible) \(d \times d\) matrices with entries in \(K\).
Note that when \(K\) is infinite \(\allmat_d(K)\) acts on any polynomial representation of \(\GL_d(K)\) by extending the domain of the defining polynomials.

Let \(d' \leq d\) and let \(\epsilon = \twobytwosmallmatrix{I_{d'}}{0}{0}{0}\in \allmat_d(K)\), a block matrix, where \(I_{d'}\) is the \(d' \times d'\) identity matrix.
Note that \(\epsilon\) is an idempotent and that \(\epsilon K\GL_d(K) \epsilon \iso K\GL_{d'}(K)\) as algebras.

\begin{definition}[Dimension reduction functor]
Suppose \(K\) is infinite  and \(d' \leq d\).
The \emph{dimension reduction functor} from \(d\) to \(d'\) is the functor from the category of polynomial left \(K\GL_d(K)\)-modules of degree \(n\) to the category of polynomial left \(K\GL_{d'}(K)\)-modules of degree \(n\) defined by left multiplication by \(\epsilon\).
\end{definition}

\begin{proposition}
\label{prop:properties_of_dimension_change}
Suppose \(K\) is infinite and \(d' \leq d\).
\begin{enumerate}[(i)]
    \item\label{item:dimension_change_exact}
The dimension reduction functor is exact.
    \item\label{item:dimension_change_on_nabla}
For any \(K\GL_d(K)\)-module \(V\), we have \(\epsilon \nabla^\lambda(V) \iso \nabla^{\lambda}(\epsilon V)\).
    \item\label{item:dimension_change_on_gtensor}
For any \(KS_n\)-module \(U\), we have \(\epsilon \gtensor^{d} (U) \iso \gtensor^{d'}(U)\).
\end{enumerate}
\end{proposition}

\begin{proof}
For part \ref{item:dimension_change_exact}, see \cite[(6.2a)]{green2006polynomial}.
Part \ref{item:dimension_change_on_nabla} is clear, and in the case of \(V=E\) is noted in \cite[Remark following (6.5f)]{green2006polynomial}.
For part \ref{item:dimension_change_on_gtensor}, let \(E'\) denote the natural \(K\GL_{d'}(K)\)-module, and observe that \(\epsilon E \iso E'\) and that furthermore \(\epsilon (E^{\tensor n}) \iso (E')^{\tensor n}\).
The proposition then follows by the definition of \(\gtensor\).
\end{proof}

This proposition tells us that, informally, the structure of \(\gtensor^{d}(S^\lambda)\) is independent of \(d\).
More precisely, we have the following corollary.

\begin{corollary}
\label{cor:dimension_independent_structure}
Suppose \(K\) is infinite.
Let \(\underline{\mu} = (\mu^{(1)}, \ldots, \mu^{(r)})\) be the sequence of labels for the simple modules in a composition series of \(\gtensor^d(S^\lambda)\) for some fixed \(d \geq n\).
Then \(\underline{\mu}\) is also the sequence of labels for the simple modules in a composition series for \(\gtensor^{d'}(S^\lambda)\) for any \(d'\) (after excluding the labels for zero modules).
\end{corollary}

\section{Quotient construction of the image of the Specht module}
\label{section:quotient_construction}

In this section we present an explicit model for \(\gtensor(S^\lambda)\) in all characteristics.
The isomorphism \(\gtensor(S^\lambda) \iso \dWeyl{\lambda}\) stated in \Cref{thm:mainresultnot2} for characteristics not \(2\) follows immediately.

\subsection{Skew column tabloids}
\label{subsection:skew_column_tabloids}

To describe the image of the Specht module under \(\gtensor\), we require a modified notion of column tabloid, which we introduce in this subsection.

Recall that in \Cref{subsection:tableaux_and_tabloids} we defined a column tabloid as an element \(t + \Jalt\) in the quotient \(\Tbx{\lambda}{V}/\Jalt\), where \(\Jalt\) is the subspace
\[
    \Jalt = \spangle{x \in \Tbx{\lambda}{V}}{\text{\(x \ppa \tau = x\) for some transposition \(\tau \in \CPP(\lambda)\)}}{K}.
\]
Consider instead the quotient by the subspace
\[
    \Jsk = \spangle{x \ppa \tau - x\sign{\tau}}{\text{\(x \in \Tbx{\lambda}{V}\), \(\tau \in \CPP(\lambda)\) is a transposition}}{K}.
\]
Note that \(\Jsk \subseteq \Jalt\) with equality if \(\charac{K} \neq 2\).
In characteristic \(2\), the additional elements of \(\Jalt\) are the tableaux with repeated entries in a column; that is:
\[\Jalt = \Jsk + \spangle{t \in \Tbx{\lambda}{V}}{\text{\(t\) has a repeated entry in a column}}{K}.\]

\begin{definition}[Skew column tabloid]
\label{def:skew_column_tabloid}
The \emph{skew column tabloid} corresponding to a tableau \(t\) is the element \(t + \Jsk\) in the quotient \(\Tbx{\lambda}{V}/\Jsk\). 
\end{definition}


We write the skew column tabloid corresponding to a tableau \(t\) as \(\sct{t}\), and draw a skew column tabloid by deleting the horizontal lines from a drawing of the corresponding tableau and double-striking the vertical lines, as depicted below in the case \(\lambda = (3,2)\).
\[
    t = \ytab{124,35} 
        \quad \implies\quad
    \sct{t} = \yskewcoltab{124,35}
\]

Depending on the characteristic, the space of skew column tabloids is isomorphic as a \(KG\)-module either to an exterior power or symmetric power:
\[
    \faktor{\Tbx{\lambda}{V}}{\Jsk} \iso \begin{cases}
    \Wedge^{\lambda'} V & \text{if \(\charac{K} \neq 2\),} \\
    \Sym^{\lambda'} V & \text{if \(\charac{K} = 2\).}
    \end{cases}
\]
For convenience, we define the \emph{skew symmetric power} \(\Sk^{\blank}\) to be the symmetric power \(\Sym^{\blank}\) in characteristic \(2\) and the exterior power \(\Wedge^{\blank}\) otherwise.
We use \(\skewtabs{\lambda}{V}\) to denote the space of skew column tabloids.

As is already clear, the definitions of alternating column tabloids and skew column tabloids agree in characteristics other than \(2\).
The definitions also agree if we restrict to tableaux of symmetric type.

Both alternating and skew column tabloids have the property that, for \(\sigma \in \CPP(\lambda)\), the equalities
\begin{align*}
    \act{t \cdot \sigma} &= \act{t} \sign{\sigma} \\
    \sct{t \cdot \sigma} &= \sct{t} \sign{\sigma}
\end{align*}
hold.
The key difference between the two definitions of tabloids is that alternating column tabloids furthermore have the property that if \(t\) has a repeated entry in a column then \(\act{t} = 0\), whereas skew column tabloids do not have this property in characteristic \(2\).
It is for these properties that the tabloids are named: an alternating column tabloid resembles an alternating multilinear form, whereas  a skew column tabloid resembles a skew symmetric multilinear form.

There is a surjection \(\qmap \colon \skewtabs{\lambda}{V} \to \alttabs{\lambda}{V}\) defined by \(K\)-linear extension of
\[
    \qmap \colon \sct{t} \mapsto \act{t}.
\]
This map is easily seen to be \(G\)-equivariant.
The kernel of \(\qmap\) is the subspace spanned by skew column tabloids with repeated column entries.

We define skew Garnir relations analogously to Garnir relations, as follows.

\begin{definition}[Skew Garnir relations]
\label{def:skew_Garnir_relation}
Let \((t, A, B)\) and \(\mathcal{S}\) be as in the definition of a Garnir relation (\Cref{def:Garnir_relation}): \(t\) is a tableau of shape \(\lambda\) with entries in \(\entryset\), \(1 \leq j < j' \leq \lambda_1\), \(A \subseteq \col_j(\lambda)\) and \(B \subseteq \col_{j'}(\lambda)\) are such that \(\abs{A} + \abs{B} > \lambda'_j\), and \(\mathcal{S}\) is a set of left coset representatives for \(S_{A} \times S_B\) in \(S_{A \sqcup B}\).
Then the \emph{skew Garnir relation} labelled by \((t, A, B)\) is
\[
    \skGRel{(t,A,B)} = \sum_{\tau \in \mathcal{S}} \sct{t \ppa \tau} \sign \tau.
\]
Let \(\skGR{\lambda}{V}\) denote the subspace 
of \(\skewtabs{\lambda}{V}\) which is spanned by the Garnir relations.

If we wish to emphasise that a Garnir relation as defined in \Cref{def:Garnir_relation} is an element of \(\alttabs{\lambda}{V}\) and not a skew Garnir relation, we describe it as an \emph{alternating Garnir relation}. 
\end{definition}

Just as for the alternating Garnir relations, a skew Garnir relation does not depend on the choice of coset representatives, and the \(K\)-subspace \(\skGR{\lambda}{V}\) is moreover a \(KG\)-submodule because the group action commutes with the place permutation action.

We likewise define certain distinguished skew Garnir relations.

\begin{definition}[Skew snake relations and basic skew snake relations]
\label{def:skew_snake_relation}
A skew Garnir relation is called a \emph{skew snake relation} under the same conditions described for Garnir relations in \Cref{def:snake_relation}: if, in the notation of \Cref{def:skew_Garnir_relation}, \(j' = j+1\) and there exists \(i\) such that \(A = \setbuild{(x,j)}{i \leq x \leq \lambda'_j}\) and \(B = \setbuild{(x,j')}{1 \leq x \leq i}\); in this case, we may also label the Garnir relation by \((t, i, j)\).
Given a function \(\Phi\) on column semistandard tableaux which are not row semistandard whose output on such a tableau \(t\) is a box \((i,j)\) such that \(t(i,j) > t(i,j+1)\), a skew snake relation labelled by \((t,i,j)\) is called \emph{\(\Phi\)-basic} if \(t\) is column semistandard but not row semistandard and \((i,j) = \Phi(t)\) (that is, under the same conditions described for Garnir relations in \Cref{def:basic_snake_relation}, with ``column standard'' replaced with ``column semistandard'').
\end{definition}

The image of a skew Garnir relation \(\skGRel{(t,A,B)}\) under \(\qmap \colon \skewtabs{\lambda}{V} \to \alttabs{\lambda}{V}\) is of course the Garnir relation \(\GRel{(t,A,B)}\).
However, \(\skGRel{(t,A,B)}\) may have nonzero summands which vanish under \(\qmap\).

\subsection{The image of the Specht module under the inverse Schur functor}
\label{subsection:iso_with_Weyl_module}

We can now use the skew column tabloids and the skew Garnir relations to model the image of the Specht module under the inverse Schur functor.

Recall that \(E\) denotes the natural representation of \(\GL_d(K)\), \(\natpermrep\) denotes the natural permutation representation of \(S_n\), and that \(S^\lambda = \nabla^\lambda_\sym W\).
We view the basis of \(W\) as \([n]\); let \(\entryset\) denote a basis for \(E\).

\begin{lemma}
\label{lemma:gtensor_applied_to_sym_ses}
Let \(n\) and \(d\) be any integers.
\begin{enumerate}[(i)]
    \item\label{item:gtensor_on_alt_tabs_of_sym_type}
There is an isomorphism \(\gtensor(\symalttabs{\lambda}{\natpermrep}) \iso \skewtabs{\lambda}{E}\).
    \item\label{item:image_of_ses}
There is a short exact sequence
\begin{center}
\begin{tikzcd}
0 \arrow[r] & \skGR{\lambda}{E} \arrow[r] & \skewtabs{\lambda}{E} \arrow[r] & \gtensor(S^\lambda) \arrow[r] & 0
\end{tikzcd}
\end{center}
in the category of \(K\GL_d(K)\)-modules.
\end{enumerate}
\end{lemma}

\begin{proof}
{[\ref{item:gtensor_on_alt_tabs_of_sym_type}]}
Given a pure tensor \(x = x_1 \tensor \cdots \tensor x_n \in E^{\tensor n}\) whose factors are basis elements in \(\entryset\), and given also a tableau \(u\) of shape \(\lambda\) with entries in \([n]\), let \(x_u\) denote the tableau of shape \(\lambda\) with entries in \(\entryset\) defined by
\[
    x_u (b) = x_{u(b)}
\]
for all \(b \in \Y{\lambda}\).

Fix any tableau \(s\) of symmetric type with entries in \([n]\).
Then there are mutually inverse \(K\GL_d(K)\)-isomorphisms \(\phi \colon \gtensor(\symalttabs{\lambda}{\natpermrep}) \to \skewtabs{\lambda}{V}\) and \(\psi \colon \skewtabs{\lambda}{V} \to \gtensor(\symalttabs{\lambda}{\natpermrep})\) given by \(K\)-linear extension of
\begin{align*}
    \phi(x \tensor_{KS_n} \act{u}) &= \sct{x_u}
\intertext{and}
    \psi(\sct{t}) &= \bigotimes_{i \in [n]} t(s^{-1}(i)) \tensor_{KS_n} \act{s}
\end{align*}
for all elements \(x\) and \(u\) as above and all tableaux \(t\) with entries in \(\entryset\).
For example, with \(\lambda = (3,2)\) there is a correspondence between elements
\[
x_1 \tensor \cdots \tensor x_5 \tensor_{KS_n} \ycoltab{124,35}
\ \;\leftrightarrow\;
\yskewcoltab{{x_1}{x_2}{x_4},{x_3}{x_5}}
\]
under \(\phi\) and \(\psi\).

The existence of \(\phi\) can be shown using the universal property of the tensor product.
Showing that \(\phi\) and \(\psi\) are well-defined and \(\GL_d(K)\)-equivariant and that they are left- and right-inverses is an exercise in bookkeeping (using the choice of tableau \(s\) to translate between \(S_{\Y{\lambda}}\) and \(S_{n}\)).

[\ref{item:image_of_ses}]
By \Cref{prop:short_exact_sequence_of_symmetric_type}, there is a short exact sequence
\begin{center}
\begin{tikzcd}
0 \arrow[r] & \symGR{\lambda}{\natpermrep} \arrow[r, "\iota"] & \symalttabs{\lambda}{\natpermrep} \arrow[r, "e|_{\sym}"] & S^\lambda \arrow[r] & 0
\end{tikzcd}
\end{center}
where \(\iota\) denotes the inclusion map.
Since \(\gtensor\) is right-exact, applying it to this sequence we obtain an exact sequence ending
\begin{center}
\begin{tikzcd}
\arrow[r] & \gtensor(\symGR{\lambda}{\natpermrep}) \arrow[r, "\gtensor(\iota)"] & \gtensor(\symalttabs{\lambda}{\natpermrep}) \arrow[r, "\gtensor(e|_{\sym})"] & \gtensor(S^\lambda) \arrow[r] & 0.
\end{tikzcd}
\end{center}
Applying the isomorphism \(\phi\) from \ref{item:gtensor_on_alt_tabs_of_sym_type}, we have a short exact sequence
\begin{center}
\begin{tikzcd}
0 \arrow[r] & \im \phi\gtensor(\iota) \arrow[r] & \skewtabs{\lambda}{V} \arrow[r] & \gtensor(S^\lambda) \arrow[r] & 0
\end{tikzcd}
\end{center}
and it suffices to show that \(\im \phi\gtensor(\iota) = \skGR{\lambda}{V}\).

The image \(\im\phi\gtensor(\iota)\) is spanned by elements of the form
\[
    \phi(x \tensor_{KS_n} \mathsf{R}_{(t,A,B)})
\]
where \(t\) is a tableau with entries in \([n]\), \(A\) and \(B\) are subsets of \(\Y{\lambda}\) as in the definition of a Garnir relation, and \(x\) is a pure tensor whose factors are basis elements of \(E\).
Fix such \(t\), \(A\), \(B\) and \(x\), and let \(\mathcal{S}\) be a set of left coset representatives for \(S_{A} \times S_{B}\) in \(S_{A \sqcup B}\).
Then, using that \(x_{t \ppa \sigma} = x_t \ppa \sigma\) for any \(\sigma \in S_{\Y{\lambda}}\), we have
\begin{align*}
\phi(x \tensor_{KS_n} \mathsf{R}_{(t,A,B)})
    &= \sum_{\tau \in \mathcal{S}} \phi(x \tensor_{KS_n} \act{t \ppa \tau} \sign\tau) \\
    &= \sum_{\tau \in \mathcal{S}} \sct{x_{t \ppa \tau}} \sign\tau \\
    &= \sum_{\tau \in \mathcal{S}} \sct{x_t \ppa \tau} \sign\tau,
\end{align*}
which is a skew Garnir relation labelled by \((x_t,A,B)\).
Since also any tableau with entries in \(\entryset\) can be written in the form \(x_t\) for suitable \(x\) and \(t\), we have that \(\im \phi\gtensor(\iota) = \skGR{\lambda}{V}\) as required.
\end{proof}

\begin{proposition}
\label{prop:commutative_diagram_for_gtensor}
The following diagram in the category of \(K\GL_d(K)\)-modules is commutative with exact rows and exact columns.
In particular, there is a surjection \(\gtensor(S^\lambda) \onto \dWeyl{\lambda}\) which is an isomorphism if and only if \(\ker q \subseteq \skGR{\lambda}{E}\).
\begin{center}
\begin{tikzcd}
        & 0 \arrow[d]
        & 0 \arrow[d]
        & 0 \arrow[d]
        & \\
    0 \arrow[r]
        & \ker q|_{\Gar} \arrow[r, hook] \arrow[d, hook] 
        & \skGR{\lambda}{E} \arrow[r, "q|_{\Gar}"] \arrow[d, hook]
        & \GR{\lambda}{E} \arrow[r] \arrow[d, hook]
        & 0 \\
    0 \arrow[r] 
        & \ker q \arrow[r, hook] \arrow[d]
        & \skewtabs{\lambda}{E} \arrow[r, "q"] \arrow[d]
        & \alttabs{\lambda}{E} \arrow[r] \arrow[d, "e"]
        & 0 \\
     0 \arrow[r]
        & \faktor{\ker q}{\ker q|_{\Gar}} \arrow[r] \arrow[d]
        & \mathcal{G}_{\tensor}(S^\lambda) \arrow[r] \arrow[d]
        & \dWeyl{\lambda} \arrow[r] \arrow[d]
        & 0 \\
        & 0 
        & 0
        & 0
\end{tikzcd}
\end{center}
\end{proposition}

\begin{proof}
Clearly the first column and the first two rows are exact and the top two squares commute. The second and third columns are exact by \Cref{cor:short_exact_sequence_for_nabla} and \Cref{lemma:gtensor_applied_to_sym_ses}.

The existence of the maps in the third row and the commutativity of the bottom two squares follow from the universal properties of the objects in the third row as cokernels.
The map from \(\gtensor(S^\lambda)\) to \(\dWeyl{\lambda}\) is surjective by surjectivity of \(eq\) and the commutativity of the diagram.
Exactness at the remaining two objects in the third row follows from (a degenerate case of) the snake lemma.

We see immediately from the diagram that the surjection \(\gtensor(S^\lambda) \onto \dWeyl{\lambda}\) is injective if and only if \(\ker q = \ker q|_{\Gar}\), or equivalently \(\ker q \subseteq \skGR{\lambda}{E}\).
\end{proof}

From this proposition we can immediately identify the image of the Specht module in characteristics other than \(2\) (when \(\qmap\) is an isomorphism), obtaining the first main result of this paper.

\mainresultnottwo*

\section{Combinatorics of skew Garnir relations}
\label{section:skew_Garnir_relations}

The goal of this section is to identify for which partitions the necessary and sufficient condition from \Cref{prop:commutative_diagram_for_gtensor} holds in characteristic \(2\); that is, for which partitions a skew column tabloid with a repeated column entry can be written as a linear combination of skew Garnir relations.

Our classification of partitions, \Cref{prop:characterisation_of_kerq_in_skGR}, is proved in \Cref{subsection:writing_tabloids_with_Garnir_relations}.
To reach it, we must first identify a spanning set for the space of skew Garnir relations, which we do in \Cref{subsection:skew_GR_spanning_set}.

Although our application concerns representations of \(\GL_d(K)\), in this section the group action is irrelevant, so we state our results for an arbitrary \(d\)-dimensional representation \(V\) of an arbitrary group \(G\), with ordered basis \(\entryset\) viewed as \([d]\).

\subsection{Spanning set for the skew Garnir relations}
\label{subsection:skew_GR_spanning_set}

We begin by observing that the basic skew snake relations are not sufficient to span the space of skew Garnir relations.


\begin{example}
\label{eg:tabloid_equals_Garnir_relation}
Suppose \(\lambda = (2,1)\) and \(t = \ytabsmall{11,1}\). Then
\[
    \skR_{(t,1,1)} = \yskewcoltab{11,1} + \yskewcoltab{11,1} + \yskewcoltab{11,1} = \yskewcoltab{11,1}
\]
which is nonzero (in characteristic \(2\)).
However, \(\mathsf{R}_{(t,1,1)} = 0\).
\end{example}

\Cref{eg:tabloid_equals_Garnir_relation} also illustrates that the basic skew snake relations do not span the skew Garnir relations.
Indeed, the relation above is the unique skew Garnir relation for tableaux of this weight, but it is not a basic skew snake relation.

The additional skew snake relations we require to form a spanning set are defined below.
To prove that they span, we introduce additional symbols to force the tableaux to have distinct entries, then use the basis for Garnir relations of symmetric type identified in \Cref{prop:basis_for_Garnir_relations_of_sym_type} and map back down to the case of interest.

\begin{definition}[Supplementary skew snake relations]
\label{def:supplementary_skew_snake_relations}
A skew snake relation labelled by \((t,i,j)\) is called \emph{supplementary} if \(t\) is row-and-column semistandard and \(t(i,j) = t(i,j+1)\).
\end{definition}

\begin{proposition}
\label{prop:spanning_set_for_skew_Garnir_relations}
The basic and supplementary skew snake relations together span \(\skGR{\lambda}{V}\).
\end{proposition}

\newcommand{\projmap}{\pi_1}
\newcommand{\lifted}[1]{#1^\vee}
\newcommand{\hatprojmap}{\hat{\pi}_1}

\begin{proof}
Recall that the set in which our tableaux take entries is an ordered basis \(\entryset\) of \(V\).
Let \(\lifted{\entryset} = \entryset \times [n]\), ordered lexicographically, and let \(\lifted{V}\) be the \(K\)-vector space with basis \(\lifted{\entryset}\).
Let \(\projmap \colon \lifted{\entryset} \to \entryset\) be the surjection defined by \(\projmap(m,r) = m\).
Extend \(\projmap\) to a map on tableaux by acting entrywise.
This map is also surjective, and remains so on restriction to tableaux of symmetric type: given any tableau \(t\) with entries in \(\entryset\), there exists a tableau \(\lifted{t}\) with entries in \(\lifted{\entryset}\) such that \(\projmap(\lifted{t}) =t\), formed by replacing each \(m \in \entryset\) with \((m,r)\) for some \(r \in [n]\), and these \(r \in [n]\) can be chosen such that all entries in \(\lifted{t}\) are distinct.

The map \(\projmap\) induces a \(K\)-linear surjection \( \hatprojmap \colon \symGR{\lambda}{\lifted{V}} \to \skGR{\lambda}{V}\) defined by sending each column tabloid of symmetric type \(\act{\lifted{t}}\) to the skew column tabloid \(\sct{\projmap(\lifted{t})}\);
that is, \(\hatprojmap(\mathsf{R}_{(\lifted{t},A,B)}) = \skR_{(\projmap(\lifted{t}),A,B))}\) for any label \((\lifted{t},A,B)\) for a Garnir relation in \(\symGR{\lambda}{\lifted{V}}\).
This is well-defined because for tableaux of symmetric type \(\lifted{t}_1\) and \(\lifted{t}_2\), there is equality \(\act{\lifted{t}_1} = \pm \act{\lifted{t}_2}\) if and only if \(\lifted{t}_1\) and \(\lifted{t}_2\) have the same column sets (this is not the case for general tableaux: we may have equality \(\act{\lifted{t}_1} = 0 = \act{\lifted{t}_2}\) in \(\alttabs{\lambda}{\lifted{V}}\) despite an inequality \(\sct{\projmap(\lifted{t}_1)} \neq \sct{\projmap(\lifted{t}_2)}\) in \(\skewtabs{\lambda}{V}\), when the tableaux have distinct column sets but some repeated column entries).

Let \(\Phi\) be the function with respect to which we consider skew snake relations in \(\skGR{\lambda}{V}\) basic.
Choose a function \(\lifted{\Phi}\) to consider snake relations in \(\symGR{\lambda}{\lifted{V}}\) basic with respect to, chosen with the property that \(\lifted{\Phi}(\lifted{t}) = \Phi(\projmap(\lifted{t}))\) whenever \(\Phi(\projmap(\lifted{t}))\) is defined (that is, whenever \(\projmap(\lifted{t})\) is not row semistandard).
Indeed this is possible: when it is defined, the box \(\Phi(\projmap(\lifted{t})) = (i,j)\) satisfies \(\lifted{t}(i,j) > \lifted{t}(i,j+1)\) by considering the first value of the pair in each box (that is, the image under \(\projmap\)).

\Cref{prop:basis_for_Garnir_relations_of_sym_type} tells us that, in \(\symGR{\lambda}{\lifted{V}}\), the \(\lifted{\Phi}\)-basic snake relations of symmetric type form a basis.
Therefore the image of this set under \(\hatprojmap\) is a spanning set for \(\skGR{\lambda}{V}\). 
It suffices to show that this image is the union of the sets of basic and supplementary skew snake relations.

Consider a skew snake relation \(\skR_{(t,i,j)} \in \skGR{\lambda}{V}\) which is either \(\Phi\)-basic or supplementary. 
We aim to show there exists a tableau \(\lifted{t}\) with entries in \(\lifted{\entryset}\) such that \(\projmap(\lifted{t}) = t\) and \((i,j) = \lifted{\Phi}(\lifted{t})\) (and hence \(\mathsf{R}_{(\lifted{t},i,j)}\) is \(\lifted{\Phi}\)-basic and its image under \(\hatprojmap\) is \(\skR_{(t,i,j)}\)).
If \(\skR_{(t,i,j)}\) is \(\Phi\)-basic, then choose any \(\lifted{t}\) such that \(\projmap(\lifted{t}) = t\);
since \(t\) is not row semistandard, neither is \(\lifted{t}\), and so by choice of \(\lifted{\Phi}\) we have \(\lifted{\Phi}(\lifted{t}) = \Phi(t) = (i,j)\).
If \(\skR_{(t,i,j)}\) is supplementary, then \(t\) is row-and-column semistandard, and so for any choice of \(\lifted{t}\) such that \(\projmap(\lifted{t}) = t\) we have that the first values of the entries of \(\lifted{t}\) weakly increase along rows and columns.
Choose the second values of the entries of \(\lifted{t}\) such that \(\lifted{t}(i,j) > \lifted{t}(i,j+1)\) and such that elsewhere the second values strictly increase along rows and columns (for example, by filling in the entries left to right of each row in turn, then swapping the entries of \((i,j)\) and \((i,j+1)\)).
Then \(i\) and \(j\) are unique such that \(\lifted{t}(i,j) > \lifted{t}(i,j+1)\), and hence \(\lifted{\Phi}(\lifted{t}) = (i,j)\).

Now consider a \(\lifted{\Phi}\)-basic snake relation \(\mathsf{R}_{(\lifted{t},i,j)} \in \symGR{\lambda}{\lifted{V}}\).
We aim to show that the skew snake relation \(\hatprojmap(\mathsf{R}_{(\lifted{t},i,j)}) = \skR_{(\projmap(\lifted{t}),i,j)}\) is either \(\Phi\)-basic or supplementary.
If \(\projmap(\lifted{t})(i,j) > \projmap(\lifted{t})(i,j+1)\), then \(\projmap(\lifted{t})\) is not row semistandard and, by choice of \(\lifted{\Phi}\), we have that \((i,j) = \lifted{\Phi}(\lifted{t}) = \Phi(\projmap(\lifted{t}))\) and hence \((\projmap(\lifted{t}),i,j)\) labels a \(\Phi\)-basic skew snake relation.
If \(\projmap(\lifted{t})(i,j) = \projmap(\lifted{t})(i,j+1)\), then \(\projmap(\lifted{t})\) is row semistandard (or else \(\Phi(\projmap(\lifted{t}))\) would be defined and not equal to \((i,j) = \lifted{\Phi}(\lifted{t})\), a contradiction), and so \((\projmap(\lifted{t}),i,j)\) labels a supplementary skew snake relation.
\end{proof}

The spanning set identified in \Cref{prop:spanning_set_for_skew_Garnir_relations} is in general not a basis: the supplementary skew snake relations are in general not linearly independent.
Indeed, a supplementary skew snake relation may even be zero, as evidenced in the following example.

\begin{example}\label{eg:zero_supplementary_skew_relation}
Suppose \(\lambda = (2,1,1)\) and \(t=\ytabsmall{11,2,2}\). Then
\[
    \skR_{(t,1,1)} = \yskewcoltab{11,2,2} + \yskewcoltab{11,2,2} + \yskewcoltab{12,1,2} + \yskewcoltab{12,1,2} = 0.
\]
It is interesting to note that this is the unique skew Garnir relation labelled by \(t\), and that the unique skew Garnir relations labelled by other tableaux with this weight also vanish; thus no tabloid of this weight appears with nonzero coefficient in any skew Garnir relation.
\end{example}

Nevertheless, the following analogue for skew snake relations of \Cref{lemma:ordered_expression_for_snake_relation} (from which the linear independence of the basic alternating snake relations followed) is useful.
In particular, it shows that the basic skew snake relations are linearly independent.

\begin{lemma}
\label{lemma:ordered_expression_for_skew_snake_relation}
Let \(t\) be a column semistandard tableau, and suppose \((i,j)\) in such that \(t(i, j) \geq t(i, j+1)\).
Then
\[
    \skR_{(t,i,j)} = m_t \sct{t} + \sum_{u \colless t} m_u \sct{u}
\]
for some elements \(m_u\) in the subring of \(K\) generated by \(1\).
If \(t(i,j) > t(i,j+1)\), then \(m_t = 1\).
If \(t(i,j)=t(i,j+1)\)
and \(a\) and \(b\) are the multiplicities of \(t(i,j)\) in the sets defining the Garnir relation
\[
    A = \setbuild{(x,j)}{i \leq x \leq \lambda'_j} \text{ and } B = \setbuild{(x,j+1)}{1 \leq x \leq i}
\]
respectively, then \(m_t = \binom{a+b}{a}\).
\end{lemma}

\begin{proof}
Analogously to the proof of \Cref{lemma:ordered_expression_for_snake_relation}, we observe that
\[
    t(1,j+1) \leq \ldots \leq t(i,j+1) \leq t(i,j) \leq t(i+1,j) \leq \ldots \leq t(\lambda'_j, j)
\]
and hence that \(t \ppa \sigma \colleq t\) for any \(\sigma \in S_{A \sqcup B}\) (where \(A\) and \(B\) are the sets defining the Garnir relation as in the statement of the lemma).
If \(t(i,j) > t(i,j+1)\), then \(t \ppa \sigma \colsim t\) holds if and only if \(\sigma \in S_A \times S_B\).
If \(t(i,j) = t(i,j+1)\),
then \(t \ppa \sigma \colsim t\) holds for precisely those permutations which, modulo \(S_A \times S_B\), permute only the boxes containing \(t(i,j)\). The number of cosets of such permutations is the number of ways to choose \(a\) of the \(a+b\) copies of the repeated entry to include in the left-hand column.
\end{proof}

\subsection{Containment of \texorpdfstring{\(\ker q\)}{ker(q)} in the skew Garnir relations}
\label{subsection:writing_tabloids_with_Garnir_relations}

In this subsection, we characterise when there is containment \(\ker \qmap \subseteq \skGR{\lambda}{V}\) in characteristic \(2\).
When \(V=E\), this containment is equivalent to the existence of an isomorphism \(\gtensor(S^\lambda) \iso \dWeyl{\lambda}\) by \Cref{prop:commutative_diagram_for_gtensor}.

Recall \(\qmap \colon \skewtabs{\lambda}{V} \to \alttabs{\lambda}{V}\) is the map \(\sct{t} \mapsto \act{t}\) (defined in \Cref{subsection:skew_column_tabloids}).
The kernel of \(\qmap\) is spanned by skew column tabloids with a repeated entry in a column.
We have already seen that such a tabloid may or may not lie in the space of skew Garnir relations:
\Cref{eg:tabloid_equals_Garnir_relation} exhibited a skew column tabloid in the kernel of \(\qmap\) which is equal to a skew Garnir relation, whilst
\Cref{eg:zero_supplementary_skew_relation} exhibited a skew column tabloid in the kernel of \(\qmap\) that cannot be written as a linear combination of skew Garnir relation because all relations of that weight vanish.
We further illustrate this behaviour with the following example.

\begin{example}\label{eg:hooks_with_all_entries_equal}
Fix an element \(m \in \entryset\), and let \(t\) be the tableau whose entries are all \(m\).
Provided \(\lambda\) has at least two rows, we have \(\sct{t} \in \ker q\).
Meanwhile, \(t\) is the unique tableau of its weight, so it labels all skew Garnir relations of its weight.
All summands of such relations are equal to \(\sct{t}\), so
\[
    \skR_{(t,A,B)} = \begin{cases}
        \sct{t} & \text{if the number of summands \(\binom{\abs{A} + \abs{B}}{\abs{A}}\) is odd;} \\
        0 & \text{otherwise.}
    \end{cases}
\]

Suppose \(\lambda\) is a hook partition which has at least two rows and two columns.
Let \(a \geq 2\) and \(l \geq 2\) be such that \(\lambda = (a, 1^{l-1})\).
Clearly the skew Garnir relations involving only columns of length 1 have exactly two summands and hence are zero.
The number of summands in a skew Garnir relation involving the first column is \(\binom{l+1}{1} = l+1\), which is odd if and only if \(l\) is even.
Thus \(\sct{t} \in \skGR{\lambda}{V}\) holds if and only if \(l\) is even.
\end{example}

We now proceed with classifying when \(\ker q \subseteq \skGR{\lambda}{V}\).
Recall we say that a partition is \(2\)-regular if it has no repeated (positive) parts, and that it is \(2\)-singular otherwise.

\begin{lemma}
\label{lemma:tabloids_are_snakes_when_2-regular}
Suppose \(\lambda\) is \(2\)-regular, or \(\lambda_1 = \lambda_2 \geq \lambda_3 + 2\) and \(\lambda\) minus its first part is \(2\)-regular.
Then \(\ker q \subseteq \skGR{\lambda}{V}\).
\end{lemma}

\begin{proof}
First note that \(\ker q = 0\) if \(\lambda\) has exactly one row;
the lemma holds trivially in this case, so we may assume that \(\lambda\) has at least two rows, and hence that there exist tableaux \(t\) such that \(\sct{t} \in \ker q\).
Since \(\lambda\) has at least two rows, if \(\lambda\) is \(2\)-regular then \(\lambda_1 > \lambda_2 > \lambda_3\); thus under either hypothesis we have \(\lambda_1 \geq \lambda_3 +2\).

Let \(t\) be a tableau such that \(\sct{t} \in \ker q\).
Then \(t\) has at least one column with repeated entries; let \(j\) be the index of the rightmost column in which \(t\) has repeated entries.
Let \(a_1\) and \(a_2\) be boxes in column \(j\) such that \(t(a_1) = t(a_2)\).
We proceed by downward induction on \(j\).

Suppose \(j > \lambda_3\).
Since \(\lambda_1 \geq \lambda_3 + 2\), there exists some \(j' \neq j\) such that \(\lambda_3 < j' \leq \lambda_1\).
We have \(\lambda' _j \leq 2\) and \(\lambda'_{j'} \leq 2\). 
Let \(b\) be any box in column \(j'\), and set \(A = \set{a_1, a_2}\) and \(B = \set{b}\) (or vice versa if \(j' < j\)).
Then \((t,A,B)\) labels a Garnir relation, and
\[
    \skR_{(t,A,B)} = \sct{t} + \sct{t \ppa (a_1\ b)} + \sct{t \ppa (a_2\ b)} = \sct{t}.
\]
Thus \(\sct{t} \in \skGR{\lambda}{V}\) as required.

Now suppose \(j \leq \lambda_3\). Since \(\lambda\) minus its first part is \(2\)-regular, we have that column \(j\) is at most one box longer than column \(j+1\).
Set \(A = \set{a_1, a_2}\) and \(B = \col_{j+1}[\lambda]\).
Then \((t,A,B)\) labels a Garnir relation, and
\[
    \skR_{(t,A,B)} = \sct{t} + \sum_{\set{b_1, b_2} \subseteq B} \sct{t \ppa (a_1\ b_1)(a_2\ b_2)}
\]
because the summands corresponding to permutations where only one box of \(A\) is moved cancel out.
The tableaux \(t \ppa (a_1\ b_1)(a_2\ b_2)\) in the above sum have a repeated entry in column \(j+1\), so by the inductive hypothesis their skew column tabloids lie in \(\skGR{\lambda}{V}\).
Hence so does \(\sct{t}\).
\end{proof}

\begin{lemma}
\label{lemma:kerq_not_contained_in_skGR_when_repeated_parts}
Suppose \(\lambda\) is such that \(\lambda\) minus its first part is \(2\)-singular, and suppose \(\abs{\entryset} \geq n-2\).
Then \(\ker q \not\subseteq \skGR{\lambda}{V}\).
\end{lemma}

\begin{proof}
Let \(N = \abs{\entryset}\), and view \(\entryset \iso [N]\).
Pick any \(r > 1\) such that \(\lambda_r = \lambda_{r+1} > 0 \).
Set \(m = 1+ \sum_{a=1}^{r-2} \lambda_a\), and let \(\alpha\) be the weight in which \(m\) has multiplicity \(\lambda_{r-1} + \lambda_r + \lambda_{r+1}\), and all other positive integers up to and including \(n+1 - (\lambda_{r-1} + \lambda_r + \lambda_{r+1})\) have multiplicity \(1\) (and all other integers have multiplicity \(0\)).
Let \(t\) be the \(\colless\)-greatest row-and-column semistandard tableau with weight \(\alpha\); this indeed exists because the required inequality \(\abs{\entryset} \geq n+1 - (\lambda_{r-1} + \lambda_{r} + \lambda_{r+1})\) follows from the assumption \(\abs{\entryset} \geq n- 2\).
Explicitly, \(t\) is defined by
\[
    t(i,j) = \begin{cases}
        j + \sum_{a=1}^{i-1} \lambda_a & \text{if \(1 \leq i \leq r-2\);} \\
        m & \text{if \(r-1 \leq i \leq r+1\);} \\
        j + m + \sum_{a=r+2}^{i-1} \lambda_a & \text{if \(r+2 \leq i \leq \lambda'_1\).}
    \end{cases}
\]
For example, if \(\lambda = (6,6,3,3,2,1)\) and \(r=3\), then \(m = 7\) and
\[\ytableausetup{nosmalltableaux}
    t = \ytab{123456,777777,777,777,89,{10}}.
\]

We aim to prove that \(\sct{t} \in \ker q \setminus \skGR{\lambda}{V}\). Clearly \(\sct{t} \in \ker q\).
To show that \(\sct{t}\) is not an element of \(\skGR{\lambda}{V}\), we require the following property of \(t\).
Given a skew Garnir relation, we say the \emph{leading tableau} of the relation is the \(\colless\)-greatest column semistandard tableau whose tabloid has nonzero coefficient in the relation.

\begin{subclaim}
\label{subclaim:t_greater_than_all_supplementary_relations}
If \(u\) is the leading tableau of a supplementary skew snake relation and is of weight \(\alpha\), then \(u \colless t\) (where \(\alpha\) and \(t\) are as defined above).
\end{subclaim}
\begin{proof}
Consider a supplementary skew snake relation labelled by \((s,i,j)\) (so that in particular \(s\) is row-and-column semistandard).
The leading tableau of this skew snake relation is at most \(s\) by \Cref{lemma:ordered_expression_for_skew_snake_relation}.
If \(s\) is of weight \(\alpha\), then by maximality of \(t\) we have \(s \colless t\) or \(s=t\).
Thus it remains only to show that \(t\) is not  the leading tableau of a supplementary skew snake relation labelled by \((t,i,j)\).

Consider the sets \(A\) and \(B\) defining the skew Garnir relation \(\skR_{(t,i,j)}\).
Let \(a\) be the multiplicity of \(t(i,j)\) in \(A\) and \(b\) the multiplicity of \(t(i,j)\) in \(B\).
Using \Cref{lemma:ordered_expression_for_skew_snake_relation}, the coefficient of \(\sct{t}\) in \(\skR_{(t,i,j)}\) is \(\binom{a+b}{a}\), so we are required to show that \(\binom{a+b}{a}\) is even.

By construction of \(t\), a supplementary skew snake relation labelled by \((t,i,j)\) has \(i \in \set{r-1, r, r+1}\).
We assess each possibility:
\begin{itemize}
    \item
if \(i=r-1\) and \(j \leq \lambda_r\), then \(a=3\) and \(b=1\), and indeed \(\binom{4}{3} = 4\) is even;
    \item
if \(i=r-1\) and \(j > \lambda_r\), then \(a=b=1\), and indeed \(\binom{2}{1} = 2\) is even;
    \item
if \(i=r\), then \(a=b=2\), and indeed \(\binom{4}{2} = 6\) is even;
    \item
if \(i=r+1\), then \(a=1\) and \(b=3\), and indeed \(\binom{4}{1} = 4\) is even. \qedhere
\end{itemize}
\end{proof}

Returning to the proof of the lemma, suppose towards a contradiction that \(\sct{t} \in \skGR{\lambda}{V}\).
Then there exists some linear combination \(\gamma\) of (nonzero) basic and supplementary skew snake relations of weight \(\alpha\) such that \(\gamma = \sct{t}\).
Consider the basic and supplementary skew snake relations with nonzero coefficient in \(\gamma\), and consider the set of their (column semistandard) leading tableaux.
Let \(u\) be \(\colless\)-greatest in this set.
We cannot have \(u \colless t\) (or else \(\sct{t}\) does not occur in any of the relations with nonzero coefficient in \(\gamma\)), and so \Cref{subclaim:t_greater_than_all_supplementary_relations} says that \(u\) is not the leading tableau of a supplementary skew snake relation.
Hence \(u\) is the leading tableau of a (unique) basic skew snake relation, and furthermore labels that relation (since by \Cref{lemma:ordered_expression_for_skew_snake_relation} the leading tableau of a basic skew snake relation is its labelling tableau).

By maximality of \(u\), the basic skew snake relation labelled by \(u\) is the unique relation with nonzero coefficient in \(\gamma\) which has \(\sct{u}\) as a summand.
Thus \(\sct{u}\) has nonzero coefficient in \(\gamma\), and hence \(\sct{u} = \pm \sct{t}\).
Since \(u\) and \(t\) are both column semistandard, we have \(u = t\).
But \(t\) is row semistandard, which contradicts that \(u\) labels a basic skew snake relation.
\end{proof}

\begin{lemma}
\label{lemma:kerq_not_contained_in_skGR_when_first_two_rows_short}
Suppose \(\lambda_1 = \lambda_2 = \lambda_3 + 1\), and suppose that \(\abs{\entryset} \geq n-2\).
Then \(\ker q \not\subseteq \skGR{\lambda}{V}\).
\end{lemma}

\begin{proof}
We argue as in the proof of \Cref{lemma:kerq_not_contained_in_skGR_when_repeated_parts} with a different choice of \(t\).
Let \(N = \abs{\entryset}\), and view \(\entryset \iso [N]\).
Let \(\alpha\) be the weight in which \(1\) has multiplicity \(\lambda_1+\lambda_2+\lambda_3\), and all other positive integers up to and including \(n+1 - (\lambda_1 +\lambda_2 + \lambda_3)\) have multiplicity \(1\) (and all other integers have multiplicity \(0\)).
Let \(t\) be the \(\colless\)-greatest row-and-column semistandard tableau with weight \(\alpha\); indeed such tableaux exist as the required inequality \(\abs{\entryset} \geq n+1 - (\lambda_{1} + \lambda_{2} + \lambda_{3})\) follows from the assumption \(\abs{\entryset} \geq n- 2\).
Explicitly, \(t\) is defined by
\[
    t(i,j) = \begin{cases}
        1 & \text{if \(1 \leq i \leq 3\);} \\
        j + 1 + \sum_{a=3}^{i-1} \lambda_a & \text{if \(4 \leq i \leq \lambda'_1\).}
    \end{cases}
\]
For example, if \(\lambda = (5,5,4,3,1)\), then
\[
    t = \ytab{11111,11111,1111,234,5}.
\]

We deduce, using \Cref{lemma:ordered_expression_for_skew_snake_relation} as in the proof of \Cref{subclaim:t_greater_than_all_supplementary_relations}, that \(t\) satisfies \(u \colless t\) for any \(u\) of weight \(\alpha\) which is the leading tableau of a supplementary skew snake relation.
Then arguing as in the final paragraphs of \Cref{lemma:kerq_not_contained_in_skGR_when_repeated_parts}, we conclude that \(\sct{t} \in \ker q \setminus \skGR{\lambda}{V}\).
\end{proof}

Combining \Cref{lemma:tabloids_are_snakes_when_2-regular,lemma:kerq_not_contained_in_skGR_when_repeated_parts,lemma:kerq_not_contained_in_skGR_when_first_two_rows_short}, we have the following characterisation of when \(\ker q \subseteq \skGR{\lambda}{V}\).

\begin{proposition}
\label{prop:characterisation_of_kerq_in_skGR}
There is containment \(\ker q \subseteq \skGR{\lambda}{V}\) if \(\lambda\) is \(2\)-regular, or if \(\lambda_1 = \lambda_2 \geq \lambda_3 + 2\) and \(\lambda\) minus its first part is \(2\)-regular.
Supposing \(\abs{\entryset} \geq n-2\), if \(\lambda\) is not of this form then \(\ker q \not\subseteq \skGR{\lambda}{V}\).
\end{proposition}

\begin{remark}
\label{remark:weaker_restriction_on_entryset_size}
In \Cref{lemma:kerq_not_contained_in_skGR_when_repeated_parts,lemma:kerq_not_contained_in_skGR_when_first_two_rows_short} and \Cref{prop:characterisation_of_kerq_in_skGR}, the restriction on \(\abs{\entryset}\) is required to ensure that we can choose a tableau with entries all distinct except for in three specified rows.
The restriction on \(\abs{\entryset}\) can be weakened if we permit dependence on \(\lambda\): in \Cref{lemma:kerq_not_contained_in_skGR_when_repeated_parts}, it is sufficient to require that \(\abs{\entryset} \geq n+1 - (\lambda_{r-1} + \lambda_r + \lambda_{r+1})\) where \(r > 1\) is minimal such that \(\lambda_r = \lambda_{r+1}\); in \Cref{lemma:kerq_not_contained_in_skGR_when_first_two_rows_short} it is sufficient to require that \(\abs{\entryset} \geq n+1 - (\lambda_{1} + \lambda_2 + \lambda_{3})\).
\end{remark}

\section{Image of the Specht module in characteristic 2}
\label{section:image_of_Specht_in_characteristic_2}

The second main result of this paper, restated below, is now clear by combining \Cref{prop:commutative_diagram_for_gtensor,prop:characterisation_of_kerq_in_skGR}.

\mainresulttwo*

The remainder of this paper establishes a variety of other results concerning \(\gtensor(S^\lambda)\) in characteristic \(2\).
The following subsections are logically independent.

In \Cref{subsection:indecomposable_Specht_modules} we use our new knowledge of the module \(\gtensor(S^\lambda)\) to deduce the indecomposability of some Specht modules in characteristic \(2\).
In \Cref{subsection:requirement_on_d}  we show that a lower bound on \(d\) that grows with \(n\) in \Cref{thm:mainresult2} is necessary.
In \Cref{subsection:subquotients} we record some restrictions on the kernel of the surjection in \Cref{thm:mainresult2} and show
that \(\gtensor(S^\lambda)\) need not have a filtration by dual Weyl modules.
In \Cref{subsection:dim_growth_of_simple_modules} we bound the dimension growth of the kernel of the surjection.
In \Cref{subsection:descriptions_in_particular_cases} we describe \(\gtensor(S^\lambda)\) for some particular partitions and values of \(d\).

\subsection{Some indecomposable Specht modules}
\label{subsection:indecomposable_Specht_modules}

We deduce the corollary to \Cref{thm:mainresult2} mentioned in the introduction identifying some indecomposable Specht modules in characteristic \(2\).

\indecSpechts*

\begin{proof}
Suppose, towards a contradiction, that \(S^\lambda\) is decomposable; write \(S^\lambda = V_1 \oplus V_2\) for \(V_1, V_2\) nonzero submodules.
Choose some \(d \geq n\).
The functor \(\gtensor\) preserves direct sums, so applying it to this decomposition and using \Cref{thm:mainresult2} we find that \(\dWeyl{\lambda} \iso \gtensor(V_1) \oplus \gtensor(V_2)\).
Note that \(\gtensor(V_1)\) and \(\gtensor(V_2)\) are nonzero, since they are mapped by \(\sfun\) to the nonzero modules \(V_1\) and \(V_2\) respectively.
This contradicts the indecomposability of \(\dWeyl{\lambda}\).
\end{proof}

\subsection{Requirement on \texorpdfstring{\(d\)}{d}}
\label{subsection:requirement_on_d}

In \Cref{thm:mainresult2}, the restriction \(d \geq n-2\) is required to ensure the existence of a certain tableau, though the restriction can be weakened if we permit dependence on \(\lambda\) (as noted in \Cref{remark:weaker_restriction_on_entryset_size}).
It is possible for the isomorphism to fail for \(d \geq n-2\) but hold for some \(d < n-2\).
Furthermore this may happen for arbitrarily large \(d\), as demonstrated by \Cref{eg:nearly_2-core_converse_fails_for_arbitrary_d} below, so a lower bound on \(d\) that grows with \(n\) is necessary.
Bearing in mind that the composition factors of these modules are independent of \(d\) (using \Cref{cor:dimension_independent_structure}), this behaviour is due to \(\gtensor(S^\lambda)\) having composition factors which are absent in \(\dWeyl{\lambda}\) but which vanish for small \(d\).

\begin{example}
\label{eg:nearly_2-core_converse_fails_for_arbitrary_d}
Suppose \(\charac{K} = 2\).
Fix \(d \in \N\).
Let \(\lambda = (d+2, d+1, d, \ldots, 2, 1, 1)\); that is, \(\lambda\) is obtained from the \(2\)-core partition of length \(d+2\) by adding a box to the first column.
Clearly \(\lambda\) minus its first part is not \(2\)-regular, so by \Cref{thm:mainresult2} we have that \(\gtensor^{d'}(S^\lambda) \not\iso \dWeyl{\lambda}\) when \(d' \geq n-2\).
However, \(\lambda\) minus its first column is \(2\)-regular, and we claim that \(\gtensor^d(S^\lambda) \iso \dWeyl{\lambda}\).

It suffices to show that \(\ker q \subseteq \skGR{\lambda}{E}\) in this case.
Let \(t\) be a tableau with entries in \(\entryset\) such that \(\sct{t} \in \ker q\).
Then \(t\) has a repeated entry in some column, and moreover must have a repeated entry in the second column: there are \(d+1\) boxes in the second column, so there are insufficiently many basis elements of \(E\) for all of them to have distinct entries.
Then the argument of \Cref{lemma:tabloids_are_snakes_when_2-regular} can be applied: we induct downward on the index of the rightmost column in which \(t\) has a repeated entry; since this index is always at least \(2\), we do not require any constraint on the first column; all other columns satisfy the condition of being at most \(1\) longer than the next, so the argument goes through.
\end{example}

\subsection{Restrictions on filtrations}
\label{subsection:subquotients}

Let \(U^\lambda\) denote the kernel of the surjection \(\gtensor(S^\lambda) \onto \dWeyl{\lambda}\), isomorphic to \(\faktor{\ker \qmap}{\ker \qmap|_{\Gar}}\).
In this subsection we record some restrictions on the structure of \(U^\lambda\) and of \(\gtensor(S^\lambda)\).
In particular, when \(d \geq n\), we find that \(U^\lambda\) does not have \(2\)-restricted composition factors and does not have a dual Weyl module as a subquotient.
We also give an example to show that \(\gtensor(S^\lambda)\) need not have a \(\nabla\)-filtration (that is, a filtration by dual Weyl modules).

\begin{proposition}
\label{prop:immediate_properties_of_U}
Suppose \(K\) is infinite and \(d \geq n\).
\begin{enumerate}[(i)]
    \item\label{item:sfun_on_U}
\(\sfun(U^\lambda) = 0\).
    \item\label{item:comp_factors_of_U_are_not_2-restricted}
If \(L^\mu(E)\) is a composition factor of \(U^\lambda\), then \(\mu\) is not \(2\)-restricted.
    \item\label{item:nablas_are_not_subquotients_of_U}
\(\nabla^\mu(E)\) is not a subquotient of \(U^\lambda\) for any partition \(\mu \vdash n\).
\end{enumerate}
\end{proposition}

\begin{proof}
Applying the exact functor \(\sfun\) to the third row of the diagram in \Cref{prop:commutative_diagram_for_gtensor}, we have a short exact sequence
\begin{center}
\begin{tikzcd}
0 \arrow[r] & \sfun(U^\lambda) \arrow[r] & \sfun\gtensor(S^\lambda) \arrow[r] & \sfun(\dWeyl{\lambda}) \arrow[r] & 0.
\end{tikzcd}
\end{center}
But \(\sfun\gtensor(S^\lambda) \iso S^\lambda \iso \sfun(\dWeyl{\lambda})\), so \ref{item:sfun_on_U} follows.
It is known that \(\sfun(L^\mu(E)) = 0\) if and only if \(\mu\) is not \(2\)-restricted \cite[(6.4a),(6.4b)]{green2006polynomial},
so \ref{item:comp_factors_of_U_are_not_2-restricted} follows from \ref{item:sfun_on_U}. 

Every dual Weyl module in characteristic \(p\) has a composition factor \(L^\mu(E)\) with \(\mu\) a \(p\)-restricted partition (this can be deduced by interpreting, as in \cite{james1980decomposition}, the decomposition matrix for \(S_n\) as a submatrix of the decomposition matrix for \(\GL_d(K)\)).
By \ref{item:comp_factors_of_U_are_not_2-restricted}, such a composition factor cannot occur in \(U^\lambda\), so \ref{item:nablas_are_not_subquotients_of_U} follows.
\end{proof}

The following example demonstrates that \(\gtensor(S^\lambda)\) need not have a \(\nabla\)-filtration.
Our strategy is to use a dimension counting argument to deduce the composition factors of \(\gtensor(S^\lambda)\) in a particular case, and then observe that this multiset of composition factors does not permit a \(\nabla\)-filtration.
The same strategy can be used to identify the composition factors of \(\gtensor(S^\lambda)\) in other small cases; the results for partitions of \(n \leq 5\) are recorded in \Cref{eg:comp_factors_for_partitions_of_n_less_than_5}.

\begin{example}[\(\gtensor(S^\lambda)\) need not have a \(\nabla\)-filtration]
\label{eg:no_nabla_filtration}
Let \(n=5\) and \(\lambda = (2,2,1)\).
Suppose \(K\) is infinite, \(\charac{K} = 2\) and \(d \geq n-1\).
We view \(\entryset \iso [d]\).
It can be shown directly that for any tableau \(t\) whose skew column tabloid lies in \(\ker q\), given any other tableau \(t'\) of the same weight there exist skew Garnir relations \(\gamma\) also lying in \(\ker q\) such that \(\sct{t} + \gamma = \sct{t'}\).
For example, if \(t = \ytabsmall{12,13,1}\) and \(t = \ytabsmall{11,13,2}\), we can choose \(\gamma = \skR_{(t,1,1)}\).
Furthermore, no skew column tabloid lies in \(\skGR{\lambda}{E}\), because all the skew snake relations have an even number of summands (and so every linear combination of snake relations is either zero or has at least two distinct column tabloids with nonzero coefficients).
Therefore, in \(U^\lambda\) there is exactly one distinct element \(\sct{t} + \ker q\restrictto{\Gar}\) for each weight of tableau that permits at least one repeated column entry, and these elements are linearly independent.
The number of such weights are enumerated in \Cref{tab:enumeration_of_weights}.
This allows us to compute \(\dim U^\lambda = \frac{1}{6}d^4 + \frac{5}{6}d^2\).

\newcommand{\nudge}{10pt}

\begin{table}[ht]
    \caption{The number of weights of tableaux with entries in \([d]\) which have at least one repeated entry in a column.}
    \centering
\begin{tabular}{ccc} \toprule
dominant weight & example tabloid & number of weights\\ \midrule
\((2,1^3)\) & \(\yskewcoltabsmall{12,13,4}\) & \(4\binom{d}{4}\) \\[\nudge]
\((2^2,1)\) & \(\yskewcoltabsmall{12,13,2}\) & \(3\binom{d}{3}\) \\[\nudge]
\((3,1^2)\) & \(\yskewcoltabsmall{12,13,1}\) & \(3\binom{d}{3}\) \\[\nudge]
\((3,2)\) & \(\yskewcoltabsmall{12,12,1}\) & \(2\binom{d}{2}\) \\[\nudge]
\((4,1)\) & \(\yskewcoltabsmall{11,12,1}\) & \(2\binom{d}{2}\) \\[\nudge]
\((5)\) & \(\yskewcoltabsmall{11,11,1}\) & \(d\) \\[\nudge]
\bottomrule
\end{tabular}
    \label{tab:enumeration_of_weights}
\end{table}

The dimensions of the simple modules for \(n=5\) can be computed from the dimensions of the dual Weyl modules (found using the hook content formula \cite[Theorem 7.21.2]{stanley2001enumcomb}) and the decomposition matrix for \(\GL_d(K)\) (see \cite[Appendix]{james1980decomposition}).
These dimensions are recorded in \Cref{tab:dimensions_of_simples} below.
By \Cref{cor:dimension_independent_structure}, the partitions labelling the composition factors of \(U^\lambda\) are independent of \(d\) for \(d \geq n\).
Thus \(\dim U^\lambda = \frac{1}{6}d^4 + \frac{5}{6} d^2\) is a positive linear combination of the dimensions in this table.

\renewcommand{\nudge}{5pt}
\newcommand{\colshrink}{4pt}
\addtolength{\tabcolsep}{-\colshrink} 

\begin{table}[ht]
    \caption{The dimensions of the simple \(K\GL_d(K)\)-modules of polynomial degree \(n=5\).}
    \centering
\begin{tabular}{@{\hspace{6pt}} c @{\hspace{15pt}} ccccccccc @{\hspace{6pt}} } \toprule
\(\lambda\) & \multicolumn{9}{c @{\hspace{6pt}}}{\(\dim L^\lambda(E)\)} \\ \midrule\rule{0pt}{5pt} 
\((1^5)\)&
    \(\frac{1}{120}d^5\)& \(-\) & \(\frac{1}{12}d^4\)   & \(+\) & \(\frac{7}{24}d^3\)   & \(-\) & \(\frac{5}{12}d^2\) & \(+\) & \(\frac{1}{5} d\)   \\[\nudge]
\((2,1^3)\) &
    \(\frac{1}{30}d^5\) & \(-\) & \( \frac{1}{6}d^4\)   & \(+\) & \(\frac{1}{6}d^3\)    & \(+\) & \(\frac{1}{6}d^2\) & \(-\) & \(\frac{1}{5} d\)    \\[\nudge]
\((2^2,1)\) &
    \(\frac{1}{30}d^5\) &       &                       & \(-\) & \(\frac{1}{3}d^3\)    & \(+\) & \(\frac{1}{2}d^2\) & \(-\) & \(\frac{1}{5} d\)    \\[\nudge]
\((3,1^2)\) &
                        &       & \(\frac{1}{6}d^4\)    & \(-\) & \(\frac{1}{2}d^3\)    & \(+\) & \(\frac{1}{3}d^2\) &      &                       \\[\nudge]
\((3,2)\) &
                        &       &                       &       & \(\frac{1}{2}d^3\)    & \(-\) & \(\frac{1}{2}d^2\) &      &                  \\[\nudge]
\((4,1)\) &
                        &       &\(\frac{1}{3}d^4\)     &       &                       & \(-\) & \(\frac{1}{3}d^2\) &       &                  \\[\nudge]
\((5)\) &
                        &       &                       &        &                      &       & \(d^2\)                                           \\
\bottomrule
\end{tabular}
    \label{tab:dimensions_of_simples}
\end{table}

\addtolength{\tabcolsep}{\colshrink} 

This allows us to deduce the composition factors of \(U^\lambda\).
Considering the coefficient of \(d^4\), we first deduce that \(L^{(3,1,1)}(E)\) must be a composition factor. 
Subtracting these dimensions and considering the highest remaining powers of \(d\) in turn, we deduce that the composition factors of \(U^\lambda\) are \(L^{(3,1,1)}(E)\), \(L^{(3,2)}(E)\) and \(L^{(5)}(E)\).

These composition factors, together with those of \(\dWeyl{\lambda}\), can then be compared with the possible composition series of dual Weyl modules (found from the decomposition matrix for \(\GL_d(K)\)).
Doing so reveals that the composition factors of \(\gtensor(S^\lambda)\) cannot be partitioned into sets of composition factors for dual Weyl modules, and hence that \(\gtensor(S^\lambda)\) has no \(\nabla\)-filtration.
\end{example}

\subsection{Dimension growth of the kernel of the quotient map}
\label{subsection:dim_growth_of_simple_modules}

In this subsection, we bound the dimension growth of \(U^\lambda\) as \(d\) varies.
We find that \(U^\lambda\) grows more slowly than \(\dWeyl{\lambda}\), so informally \(\dWeyl{\lambda}\) comprises ``most'' of \(\gtensor(S^\lambda)\).

We use big-\(O\) and big-\(\Theta\) notation: given functions \(f\) and \(g\), the statement \(f(d) = O(g(d))\) means that the function \(f\) grows asymptotically at most as quickly as \(g\), whilst \(f(d) = \Theta(g(d))\) means \(f\) grows asymptotically at the same rate as \(g\).

\begin{lemma}
\label{lemma:dim_condition}
Fix \(n\) and allow \(d\) to vary.
Let \(M\) be a \(K\)-vector space with basis labelled by (a subset of) tableaux with entries in \([d]\).
Let \(U\) be a \(K\)-subspace of \(M\).
Let \(r \geq 1\), and suppose all elements of \(U\) are linear combinations of basis elements labelled by tableaux with at most \(n-r\) distinct entries.
Then \(\dim U = O(d^{n-r})\).
\end{lemma}

\begin{proof}
Consider \(R \leq M\) the \(K\)-subspace linearly spanned by basis elements labelled by tableaux with at most \(n-r\) distinct entries.
There are at most \(\binomII{d}{n-r} = \binom{d+n-r-1}{n-r} = O(d^{n-r})\) possibilities for the multiset of entries of such a tableau (where \(\binomII{a}{b}\) denotes the number of multisubsets of size \(b\) in a set of size \(a\)), and there are at most \({(n-r)}^n\) possibilities for the arrangement of a given \((n-r)\)-multiset of entries into a tableau.
Thus \(\dim R = O(d^{n-r})\).
By assumption, \(U\) is a subspace of \(R\), and so \(\dim U = O(d^{n-r})\).
\end{proof}

\begin{proposition}
\label{prop:dim_growth_of_U}
Fix \(n\) and allow \(d\) to vary.
Then \(\dim U^\lambda = O(d^{n-1})\).
\end{proposition}

\begin{proof}
Skew column tabloids in \(U^\lambda\) have a repeated entry in a column, and so have at most \(n-1\) distinct entries;
the proposition then follows from \Cref{lemma:dim_condition}.
\end{proof}

\begin{remark}
The dimensions of the dual Weyl modules are known (given by the hook content formula \cite[Theorem 7.21.2]{stanley2001enumcomb}), and in particular \(\dim \dWeyl{\mu} = \Theta(d^n)\) for all partitions \(\mu\) of \(n\).
Thus \Cref{prop:dim_growth_of_U} tells us that \(U^\lambda\) grows more slowly than any dual Weyl module, and in particular more slowly than \(\dWeyl{\lambda} \iso \gtensor(S^\lambda)/U^\lambda\).
This fact also offers an alternative proof of \Cref{prop:immediate_properties_of_U}\ref{item:nablas_are_not_subquotients_of_U} when \(d\) is sufficiently large: for large \(d\), \(U^\lambda\) is too small to have \(\dWeyl{\mu}\) as a subquotient.
\end{remark}

\subsection{Descriptions in particular cases}
\label{subsection:descriptions_in_particular_cases}

In this subsection, we describe the module \(\gtensor(S^\lambda)\) for some particular tractable examples.
In particular, we:
\begin{itemize}
    \item
fully describe \(\gtensor(S^\lambda)\) when \(\lambda\) is a column, row, or two-row partition (\Cref{prop:columns_and_rows});
    \item
fully describe \(\gtensor(S^\lambda)\) when \(d=1\) (\Cref{prop:when_d=1});
    \item
compute the dimension of \(\gtensor(S^\lambda)\) when \(d=2\) and \(\lambda\) is a hook partition, and furthermore for hook partitions of even length identify \(\gtensor(S^\lambda)\) as a tensor product of known representations (\Cref{prop:hooks_when_d=2});
    \item
list the composition factors of \(\gtensor(S^\lambda)\) when \(\lambda\) is a partition of \(n \leq 5\) (\Cref{eg:comp_factors_for_partitions_of_n_less_than_5}).
\end{itemize}

\begin{proposition}[Columns, rows and two-row partitions]
\label{prop:columns_and_rows}
\leavevmode
\begin{enumerate}[(i)]
    \item\label{item:columns}
Suppose \(\lambda = (1^n)\) is a single column.
Then \(\gtensor(S^\lambda) \iso \Sk^n(E)\).
    \item\label{item:rows}
Suppose \(\lambda = (n)\) is a single row.
Then \(\gtensor(S^\lambda) \iso \Sym^n(E)\).
    \item\label{item:two-row}
Suppose \(\lambda = (n-m,m)\) is a two-row partition and \(\lambda \neq (1,1)\).
Then \(\gtensor(S^\lambda) \iso \nabla^{(n-m,m)}(E)\).
\end{enumerate}
\end{proposition}

\begin{proof}
When \(\lambda\) is a column, we observe that \(\skGR{\lambda}{E} = 0\) and so \(\gtensor(S^\lambda) \iso \skewtabs{\lambda}{E}\).
When \(\lambda\) consists of at most two rows (and \(\lambda \neq (1,1)\)), the claim follows from \Cref{thm:mainresultnot2,thm:mainresult2} (or, in the case of a single row, can be seen clearly from the fact that the skew Garnir relations become relations exchanging the entries along the row).
\end{proof}

It is interesting that even the case of \(d=1\) is nontrivial.
When \(d=1\), the dual Weyl module is easy to describe: \(\dWeyl{\lambda} = 0\) unless \(\lambda\) is a single row, in which case \(\dWeyl{\lambda} \iso \Sym^n E \iso  E^{\tensor n}\) of dimension \(1\).
For \(\gtensor(S^\lambda)\), again \(0\) and \(E^{\tensor n}\) are the only two possibilities, but now both can occur for partitions of arbitrary length, and the dichotomy of partitions is not straightforward to describe.

To distinguish between the two possibilities, we require the following result on the parity of binomial coefficients.

\begin{definition}
We say the binary addition of integers \(a\) and \(b\) is \emph{carry-free} if, for all \(i\), the \(i\)th binary digits of \(a\) and \(b\) are not both \(1\).
\end{definition}

\begin{lemma}
\label{lemma:parity_of_binom_coeffs}
Let \(a,b,c \in \N\).
\begin{enumerate}[(i)]
    \item\label{item:a_and_b}
The binomial coefficient \(\binom{a+b}{a}\) is odd if and only if the binary addition of \(a\) and \(b\) is carry-free.
    \item\label{item:m}
There exists \(1 \leq i \leq c-1\) such that \(\binom{c}{i}\) is odd if and only if \(c\) is not a power of \(2\).
When this is the case, the minimal \(i \geq 1\) such that \(\binom{c}{i}\) is odd is the maximal power of \(2\) that divides \(c\).
\end{enumerate}
\end{lemma}

\begin{proof}
Part \ref{item:a_and_b} is a consequence of Lucas's Theorem, as given (for example) in \cite[Lemma 22.4]{gdjames1978reptheorysymgroups}.
Part \ref{item:m} follows from part \ref{item:a_and_b} by writing \(c\) in binary.
\end{proof}

\begin{proposition}
\label{prop:when_d=1}
Suppose \(d=1\) and \(\charac{K} = 2\).
Then \(\gtensor(S^\lambda) = 0\) if and only if there exists \(1 \leq j < \lambda_1\) such that:
\begin{itemize}
    \item
\(\lambda'_j + 1\) is not a power of \(2\); and
    \item
\(\lambda'_{j+1} \geq 2^\nu\), where \(\nu \geq 0\) is maximal such that \(2^\nu\) divides \(\lambda'_j + 1\).
\end{itemize}
When \(\gtensor(S^\lambda) \neq 0\), we have \(\gtensor(S^\lambda) \iso E^{\tensor n}\).
\end{proposition}

\begin{proof}
Since \(d=1\), the set \(\entryset\) is a singleton and there is a unique tableau \(t\) with entries in \(\entryset\) (having all entries the same).
We therefore have that \(\gtensor(S^\lambda) = 0\) if and only if \(\sct{t} \in \skGR{\lambda}{E}\), and \(\gtensor(S^\lambda) \iso E^{\tensor n}\) otherwise.

All place permutations leave \(\sct{t}\) unchanged, so the skew Garnir relation labelled by sets \(A\) and \(B\) is the sum of \(\abs{S_{A \sqcup B} : S_A \times S_B}\) copies of \(\sct{t}\).
That is,
\begin{align*}
    \skR_{(t,A,B)} &= \frac{(\abs{A} + \abs{B})!}{\abs{A}!\abs{B}!} \sct{t}.
\intertext{Focusing on skew snake relations, this becomes}
    \skR_{(t,i,j)} &= \binom{\lambda'_j+1}{i} \sct{t}.
\end{align*}
Thus \(\sct{t} \in \skGR{\lambda}{E}\) if and only if there exists \(j\) such that \(\binom{\lambda'_j + 1}{i}\) is odd for some \(1 \leq i \leq \lambda'_{j+1}\).
The proposition then follows from \Cref{lemma:parity_of_binom_coeffs}\ref{item:m}.
\end{proof}

This proof of \Cref{prop:when_d=1} generalises the argument for hook partitions given in \Cref{eg:hooks_with_all_entries_equal}.
The following corollary of the proposition can in fact be deduced from that example.

\begin{corollary}[Hooks when \(d=1\)]
\label{cor:hooks_when_d=1}
Let \(a,l \geq 2\).
Suppose \(d=1\), \(\charac{K} = 2\), and \(\lambda = (a,1^{l-1})\) is a hook partition.
Then
\[
    \gtensor(S^\lambda) \iso \begin{cases}
        0 & \text{if \(l\) is even,} \\
        E^{\tensor n} & \text{if \(l\) is odd.}
    \end{cases}
\]
\end{corollary}

Our next example concerns hooks when \(d=2\).
Our description includes a \emph{Frobenius twist}.
Recall that if \(K\) is a field of characteristic \(p\), then the map \(x \mapsto x^p\) is a field endomorphism called the \emph{Frobenius endomorphism}.
This yields a group endomorphism of \(\GL_d(K)\) defined by acting entrywise.
Composing this map with the representing group homomorphism of a representation \(V\) over \(K\) yields a new representation, called the \emph{Frobenius twist} of \(V\), which we denote \(\Fr(V)\).
Given an element \(v \in V\), we denote the corresponding element of \(\Fr(V)\) by \(\Fr(v)\).

\begin{proposition}[Hooks when \(d=2\)]
\label{prop:hooks_when_d=2}
Let \(a,l \geq 2\).
Suppose \(d=2\), \(\charac{K} = 2\), and \(\lambda = (a,1^{l-1})\).
\begin{enumerate}[(i)]
    \item
    \label{item:leg_even}
Suppose \(l\) is even. Then \(\gtensor(S^\lambda) \iso \Fr(\Sym^{\frac{l}{2}-1}(E)) \tensor \Sym^{a-1}(E) \tensor \det E\), of dimension \(\frac12 al\).
    \item
    \label{item:leg_odd}
Suppose \(l\) is odd. Then \(\dim \gtensor(S^\lambda) = \frac12 (a+1)(l+1)\).
\end{enumerate}
\end{proposition}

\begin{proof}
\newcommand{\esk}[1]{e^{\Sk}(#1)}
Write \(\entryset = \set{X,Y}\), with \(X < Y\), for the basis of \(E\).
Given a tableau \(t\), write \(\esk{t}\) for the image of \(\sct{t}\) in \(\gtensor(S^\lambda)\).

We consider the spanning set for the skew Garnir relations identified in \Cref{prop:spanning_set_for_skew_Garnir_relations}, with \(\Phi\) defined on column semistandard but not row semistandard tableaux by choosing the right-most box eligible box in the first row (noting there is no other row with more than one box).

Skew Garnir relations involving only columns other than the first tell us precisely that in \(\gtensor(S^\lambda)\) the entries of the first row (except the first) can be permuted freely.

The remaining elements of our spanning set we must consider are labelled by \((t, 1, 1)\) for some \(t\), where either: \(t\) is row-and-column semistandard and \(t(1,1) = t(1,2)\); or the first column of \(t\) has all entries \(Y\), \(t(1,2) = X\), and the remainder of the first row is weakly increasing.

For \(0 \leq c \leq l\) and \(0 \leq r \leq a-1\), let \(t_{c,r}\) be the (column semistandard) tableau of shape \(\lambda\) where \(X\) appears \(c\) times in the first column and \(r\) times in the remaining columns, with the \(X\)s in the first column at the top, and the \(X\)s in the first row at the left (except possibly the first column).
The tableaux identified in the previous paragraph, labelling the snake relations we are still to consider, are all of the form \(t_{c,r}\) for some \(0 \leq c \leq l\) and \(0 \leq r \leq a-1\).
Additionally, if \(t_{c,r}\) is one of the identified tableaux and \(r=0\), then also \(c=0\).
In these cases, we have:
\begin{align}
\label{eq:remaining_Garnir_relations}
    \skR_{(t_{c,r}, 1,1)} = \begin{cases}
        (c+1)\sct{t_{c,r}} + (l-c)\sct{t_{c+1,r-1}} & \text{if \(r > 0\),} \\
        (l+1) \sct{t_{0,0}} & \text{if \(c=r=0\).}
    \end{cases}
\end{align}

{[\ref{item:leg_even}]}
Suppose \(l\) is even.
Then each relation (\ref{eq:remaining_Garnir_relations}) above has an odd total number of summands, and thus is equal to a single tabloid.
If \(\sct{t_{c,r}}\) appears as a relation, which is precisely if \(c\) is even, then its image in \(\gtensor(S^\lambda)\) is zero;
if it does not, which is precisely if \(c\) is odd, then its image in \(\gtensor(S^\lambda)\) is nonzero and is linearly independent of the images of all other tabloids of that form.
Thus
\[
\setbuild{\esk{t_{c,r}}}{\text{\(0 \leq c \leq l\), \(c\) odd, \(0 \leq r \leq a-1\)}}
\]
is a basis for \(\gtensor(S^\lambda)\). The dimension follows.

Let \(\phi \colon \gtensor(S^\lambda) \to
    \Fr(\Sym^{\frac{l}{2}-1}(E)) \tensor \Sym^{a-1}(E) \tensor \det E \) be the map defined by \(K\)-linear extension of
\begin{align*}
    \phi( \esk{t_{c,r}} ) =
    \Fr(X^{\frac{c-1}{2}}Y^{\frac{l-c-1}{2}}) \tensor X^{r}Y^{a-1-r} \tensor 1
\end{align*}
for \(0 \leq c \leq l\), \(c\) odd, \(0 \leq r \leq a-1\).
Since \(\phi\) is a bijection between bases, it is a linear isomorphism.
It is easy to verify that \(\phi\) respects the action of diagonal elements of \(\GL_2(K)\): the element \(\twobytwosmallmatrix{\alpha}{0}{0}{\beta}
\in \GL_2(K)\) acts on both \(\esk{t_{c,r}}\) and its image by multiplication by \(\alpha^{c+r}\beta^{n-c-r}\).
It then suffices to show \(\phi\) respects the action of transvections.

Let \(0 \leq c \leq l\), \(c\) odd, \(0 \leq r \leq a-1\).
Let \(g= \twobytwosmallmatrix{1}{0}{\alpha}{1} \in \GL_2(K)\) for some \(\alpha \in K\); that is, \(g\) is the transvection fixing \(Y\) and acting on \(X\) as \(g X = X + \alpha Y\).
Then
\begin{align*}
g X^{r}Y^{a-1-r}
    &= \sum_{j=0}^r \binom{r}{j} \alpha^{j} X^{r-j}Y^{a-1-r+j}
\intertext{and}
g \Fr(X^{\frac{c-1}{2}}Y^{\frac{l-c-1}{2}})
    &= \Fr(\twobytwosmallmatrix{1}{0}{\alpha^2}{1} X^{\frac{c-1}{2}}Y^{\frac{l-c-1}{2}}) \\
    &= \sum_{k=0}^{\frac{c-1}{2}} \binom{\frac{c-1}{2}}{k} \alpha^{2k} \Fr(X^{\frac{c-1}{2}-k}Y^{\frac{l-c-1}{2}+k}).
\intertext{Also \(\det g = 1\), so}
g \phi(\esk{t_{c,r}})
    &= \sum_{k=0}^{\frac{c-1}{2}} \sum_{j=0}^a \binom{c}{2k}\binom{r}{j} \alpha^{2k+j} \phi(\esk{t_{c-2k,r-j}}).
\end{align*}

Meanwhile,
\begin{align*}
g \esk{t_{c,r}}
    &= \sum_{i=0}^c \sum_{j=0}^r \binom{c}{i}\binom{r}{j} \alpha^{i+j} \esk{t_{c-i,r-j}} \\
    &= \sum_{k=0}^{\frac{c-1}{2}} \sum_{j=0}^r \binom{c}{2k}\binom{r}{j} \alpha^{2k+j} \esk{t_{c-2k,r-j}}
\end{align*}
where the second equality holds because \(\esk{t_{c-i,r-j}} = 0\) when \(i\) is odd, so we can relabel via \(i = 2k\).
Equivariance is then clear provided that \(\binom{\frac{c-1}{2}}{k} \equiv \binom{c}{2k} \pmod{2}\).
Indeed this follows from the \Cref{lemma:parity_of_binom_coeffs}\ref{item:a_and_b} by noting that the binary addition of \(a\) and \(b\) is carry-free if and only if the binary addition of \(2a\) and \(2b+1\) is carry-free.
Showing that \(\phi\) respects the action of \(\twobytwosmallmatrix{1}{\alpha}{0}{1}\) is analogous, and completes the proof.

{[\ref{item:leg_odd}]}
Suppose \(l\) is odd.
Then each relation (\ref{eq:remaining_Garnir_relations}) above has an even total number of summands, and thus is either zero or the sum of two distinct tableaux.
The \(c=r=0\) relation is clearly zero.
When \(r>0\), the relation is nonzero if and only if \(c\) is even.
We thus have that \(\esk{t_{c,r}} = \esk{t_{c+1,r-1}}\) for \(c\) even and \(r>0\), and furthermore that
\[
\setbuild{\esk{t_{c,r}}}{\text{\(0 \leq c \leq l\), \(c\) even, \(0 \leq r \leq a-1\)}}
    \sqcup
\setbuild{\esk{t_{c,a-1}}}{\text{\(0 \leq c \leq l\), \(c\) odd}}
\]
is a basis for \(\gtensor(S^\lambda)\).
The dimension follows.
\end{proof}

To finish, we record the composition factors of \(\gtensor(S^\lambda)\) when \(n \leq 5\)
, in the cases where \(\gtensor(S^\lambda)\) is not isomorphic to \(\dWeyl{\lambda}\).
The composition factors of \(\dWeyl{\lambda}\) are recorded in, for example, \cite[Appendix]{james1980decomposition}, so we record only the composition factors of \(U^\lambda\) (the kernel of the surjection \(\gtensor(S^\lambda) \onto \dWeyl{\lambda}\)).
By \Cref{cor:dimension_independent_structure}, the composition factors are independent of \(d\), though some may vanish for small \(d\).

\begin{example}[Partitions of \(n \leq 5\)]
\label{eg:comp_factors_for_partitions_of_n_less_than_5}
Suppose \(K\) is infinite and \(\charac{K} = 2\).
The dimension counting argument from \Cref{eg:no_nabla_filtration} can be used to compute the composition factors of \(U^\lambda\) for all partitions of \(n \leq 5\).

For \(n \leq 3\), the only partitions for which \(\gtensor(S^\lambda) \not\iso \dWeyl{\lambda}\) are the columns \((1^2)\) and \((1^3)\).
If \(\lambda = (1^2)\), then \(U^\lambda \iso L^{(2)}(E)\); if \(\lambda = (1^3)\), then \(U^\lambda \iso L^{(3)}(E)\).

For \(n = 4\), the partitions for which \(\gtensor(S^\lambda) \not\iso \dWeyl{\lambda}\) are \((1^4)\) and \((2,1^2)\).
For \(n = 5\), the partitions for which \(\gtensor(S^\lambda) \not\iso \dWeyl{\lambda}\) are \((1^5)\), \((2,1^3)\), \((2^2,1)\) and \((3,1^2)\).
The composition factors of \(U^\lambda\) in these cases are given in \Cref{tab:comp_factors_n=4_and_5}.

\newcommand{\fivetable}{
\begin{tabular}{ccccc}\toprule
            & \((3,1^2)\)   & \((3,2)\) & \((4,1)\) & \((5)\)   \\\midrule
\((1^5)\)   & 1             & 1         &           & 1         \\
\((2,1^3)\) &               &           & 1         &           \\
\((2^2,1)\) & 1             & 1         &           & 1         \\
\((3,1^2)\) & 1             & 2         &           & 1         \\
\bottomrule
\end{tabular}
}
\newcommand{\fourtable}{
\begin{tabular}{ccccc}\toprule
            & \((2^2)\) & \((3,1)\) & \((4)\)   \\\midrule
\((1^4)\)   & 1         & 1         & 1         \\
\((2,1^2)\) & 2         & 1         & 1         \\
\bottomrule
\end{tabular}
}

\begin{table}[ht]
\centering
\captionsetup{position=below}

\begin{subtable}[t]{.49\textwidth}
\centering
\fourtable
\vphantom{\fivetable}
\caption{\(n=4\)}
\end{subtable}%
\begin{subtable}[t]{.49\textwidth}
\centering
\fivetable
\caption{\(n=5\)}
\end{subtable}

\caption{The composition factors of \(U^\lambda\) for partitions of \(4\) and \(5\). The composition factors of \(U^\lambda\) are given by the row labelled \(\lambda\); the multiplicities of the simple module \(L^\mu(E)\) by the column labelled \(\mu\).}
\label{tab:comp_factors_n=4_and_5}
\end{table}
\end{example}

\section*{Acknowledgements}

The author is grateful to David Hemmer, John Murray, Mark Wildon and an anonymous referee for their comments on earlier versions of this paper.

\bibliographystyle{alpha}
\bibliography{references}

\end{document}